\documentclass[11pt]{amsart}
\usepackage{amsmath,amssymb,latexsym, fullpage}
\usepackage{longtable}
\newtheorem{theorem}{Theorem}[section]
\newtheorem{lemma}[theorem]{Lemma}
\newtheorem{definition}[theorem]{Definition}
\newtheorem{corollary}[theorem]{Corollary}
\newtheorem{proposition}[theorem]{Proposition}

\newcommand{\R}{\mathbb{R}}

\newcommand{\pf}{\noindent {\em Proof.}\quad}
\usepackage{graphicx, amssymb, latexsym, multicol}

\usepackage[usenames, dvipsnames]{color}
\usepackage{colordvi}

\usepackage{amsmath, multicol}
\usepackage{url}
\usepackage{hyperref}

\newcommand{\st}{S_T}
\newcommand{\qt}{Q_T}

\newcommand{\ou}{\overline{u}}
\newcommand{\ov}{\overline{v}}
\newcommand{\uu}{\underline{u}}
\newcommand{\uv}{\underline{v}}

\newcommand{\ma}{m_{\alpha}}
\newcommand{\huk}{\hat{u}^k}
\newcommand{\tpl}{\tilde{\phi}^L}

\newcommand{\uk}{u^k}
\newcommand{\vk}{v^k}
\newcommand{\wk}{w^k}

\makeatletter
 
 \@addtoreset{equation}{section}
\makeatother
\begin{document}

\title{Self-similar fast-reaction limits for  reaction-diffusion systems on unbounded domains}
\author{E.C.M. Crooks}
\address{Department of Mathematics, College of Science, Swansea University, 
Swansea SA2 8PP, U.K.}

\author{D. Hilhorst}
\address{CNRS and Laboratoire de Math\' ematiques,
 Universit\'e de Paris-Sud 11, F-91405 Orsay Cedex, France}
 
\begin{abstract}
We present a unified approach to characterising fast-reaction limits of systems of either two reaction-diffusion equations, or one reaction-diffusion equation and one ordinary differential equation, on unbounded domains, motivated by models of fast chemical reactions where either one or both reactant(s) is/are mobile. For appropriate initial data, solutions of four classes of problems each converge in the fast-reaction limit $k \to \infty$ to a self-similar limit profile that has one of four forms, depending on how many components diffuse and whether the spatial domain is a half or whole line. For fixed $k$, long-time convergence to these same self-similar profiles is also established, thanks to a scaling argument of Kamin. Our results generalise earlier work of Hilhorst, van der Hout and Peletier to a much wider class of problems, and provide a quantitative description of the penetration of one substance into another in both the fast-reaction and long-time regimes.
\end{abstract}

\date{}

\maketitle

\section{Introduction}
\noindent  Systems of the form
\begin{equation}
\label{syst-proto}
\begin{array}{rlrl}\displaystyle
u_t&= d_u u_{xx} -kuv,\quad &(x,t) \in (0, \infty) \times (0, T),
\\[.1 in]
\displaystyle v_t&=  -kuv,\quad &(x,t) \in (0, \infty) \times (0, T),  \\[.1 in]
u(0, t) & = U_0,  & \mbox{for} \quad t \in (0,T), \\[.1 in]
 u(x,0)&=0,\quad v(x,0)=V_0   &\mbox{for} \quad  x\in (0, \infty),
 \end{array}
 \end{equation}
 
 \noindent  arise in modelling chemical reactions
 $$ A+B \stackrel{k}{\rightarrow} C$$
 
 \noindent taking place in a semi-infinite region, modelled for simplicity by the one-dimensional spatial domain $(0, \infty)$ with surface $x=0$. Here $u$ and $v$ represent concentrations of a mobile chemical $A$ and immobile substrate $B$ respectively, $U_0$ and $V_0$ are positive constants, and $k$ is the (positive) rate constant of the reaction. The mobile reactant $u$ is initially not present in the domain $(0, \infty)$, the concentration of $u$ outside the domain imposes a boundary condition $u=U_0$ at $x=0$, and the immobile substrate is assumed initially to have uniform concentration $V_0$
 throughout $(0, \infty)$. Examples of where such systems can arise include modelling the penetration of radio-labelled antibodies into tumourous tissue,  or of carbonic acid into porous rock.  
The fast-reaction $k \to \infty$ limit of solutions of (\ref{syst-proto})  is both physically relevant, since, for example, the attachment of antibodies to tissue can be very fast
whereas the fact that antibodies are often relatively large makes diffusion typically slow, and  mathematically useful and interesting. In \cite{hhp1}, it was established by Hilhorst, van der Hout and Peletier   that  $k-$dependent solutions $(\uk, \vk)$ of (\ref{syst-proto}) converge as $k \to \infty$ on bounded time intervals $[0,T]$ to self-similar limit profiles $(u,v)(x/\sqrt{t})$ that satisfy a free boundary problem. This free boundary has the form $x=a\sqrt{t}$ where $a$ is a positive constant,  and separates the region in which the mobile chemical $A$ is present from that where it is absent,  thus characterising the rate at which, in the limit of fast reaction, $A$ invades  the immobile substrate $B$. Such information about how one substance penetrates into another has key applications to, for example, assessment of the effectiveness of radiotherapy or  prediction of rates of carbon dioxide sequestration.

\medskip
\noindent The modelling of other physical problems can clearly give rise to systems related to, but different from, (\ref{syst-proto}), for which the fast-reaction limit and characterisation of rates at which one substance invades another are again of interest. For instance, if both a reactant $u$ and substrate $v$ are mobile, such as when carbonic acid penetrates into water instead of rock, the substrate will diffuse, introducing a term $d_v v_{xx}$ into the model, and typically satisfy a zero flux boundary condition $v_x=0$ at the surface $x=0$. Similar models but with the half-line spatial domain $(0, \infty)$ replaced by the whole line $\R$, can arise,  for example in  neutralisation reactions where $u$ is the concentration of an acid, $v$ the concentration of a base, either both mobile or one mobile and one immobile, and the two initially separated chemicals are brought together to react \cite{Trev1, Trev1a, Trev2, TSDW}. 
The form of reaction can also be much more general than in (\ref{syst-proto}), because, for instance, chemicals $A$ and $B$ may react in the form
$$mA + n B \stackrel{k}{\rightarrow} C$$
where the stoichiometric coefficients $m, n \in \R$ are positive, which gives rise to interaction terms  $-k u^mv^n$ instead of $-kuv$, or more generally, $-kF(u,v)$, with suitable hypotheses on $F$.  
Since  reactions can exhibit fractional order kinetics \cite{Trev1, Trev2}, $m$ and $n$ need not be integers,  and thus it is important to allow the interaction term to be not necessarily Lipschitz continuous on $[0, \infty) \times [0, \infty)$.  
 Initial conditions can also be more complicated than the simple piecewise constant functions in (\ref{syst-proto}).

\medskip
\noindent Here we present a unified approach to characterising the self-similar fast-reaction limits for four different classes of problem
  covering all of the physical models above, with either one or both reactants mobile and spatial domain either   the whole line $\R$ or the half-line $(0, \infty)$, and with sets of conditions encompassing a broad range of both interaction terms and  initial conditions. Our framework  includes as a special case the results of \cite{hhp1} for the prototype problem (\ref{syst-proto}), and also  the  first extension of \cite{hhp1} in \cite{hhp2} that allows more general forms of $F$ than $-kuv$. Note that the simple form of reaction and single mobile reactant in  (\ref{syst-proto}) actually enables this particular problem  to be transformed to a single parabolic equation,  whereas both \cite{hhp1} and  the general framework presented here need alternative, more widely applicable ideas. 
Additionally,   we exploit a scaling argument to apply our results on convergence to self-similar limit profiles as $ k \to \infty$  to show that for fixed $k$, solutions converge in the long-time limit $t \to \infty$ to these same self-similar limit profiles  in a certain average sense. This enables us also to provide rigorous justification for some limiting self-similar profiles derived previously by asymptotic methods by Trevelyan et al \cite{TSDW} in the context of long-time behaviour of reaction fronts in two-layer systems, and in fact, the asymptotic work of \cite{Trev1, Trev1a, Trev2, TSDW}, together with \cite{hhp1, hhp2}, was central to the motivation for our work.

  \bigskip
  \noindent We treat two pairs of problems, depending on whether the spatial domain is  $\R$ or $(0, \infty)$. 
   The first pair is defined on the strip $\qt = \{(x,t): x \in \R, 0<t<T\}$, and the system considered is
$$
(P_1^k) \left\{\begin{array}{rlrl}\displaystyle
u_t&= d_u u_{xx} -kF(u,v)\quad &\mbox{in}~\qt,
\\[.1 in]
\displaystyle v_t&= d_v v_{xx} -kF(u,v)\quad &\mbox{in}~\qt, \\[.1 in]
 u(x,0)&=u^k_0(x),\quad v(x,0)=v^k_0(x)   &\qquad\mbox{for} \quad x\in \R,
\end{array}\right.$$

\medskip
\noindent  where we define
$$ u_0^{\infty} : = \left\{ \begin{array}{ll} U_0 & \mbox{for} \; x<0,\\
0 & \mbox{for} \; x>0, \end{array} \right., \;\; 
 v_0^{\infty} : = \left\{ \begin{array}{ll} 0 & \mbox{for} \; x<0,\\
V_0 & \mbox{for} \; x>0, \end{array} \right.$$

\noindent with   $U_0$, $V_0$ positive constants, 
and choose the initial data $u^k_0, v^k_0 \in C^2(\R)$ such that  $0 \leq u^k_0 \leq M$, $0 \leq v^k_0 \leq M$ for some $M \geq \max\{ U_0, V_0 \}$, 
 $$u^k_0 (x) \rightarrow U_0, 0 \;\; \mbox{and} \;\; v^k_0 (x) \rightarrow 0, V_0 \;\; \mbox{as} \;\; x \rightarrow - \infty, \infty \;\; \mbox{resp.},$$
$$\| u^k_0 - u_0^{\infty} \|_{L^1(\R)} < \infty,  \;\; \| v^k_0 - v_0^{\infty} \|_{L^1(\R)} < \infty,$$
$$k \mapsto u_0^k - u_0^{\infty}, \;\; k \mapsto v_0^k - v_0^{\infty}\;\; \mbox{belong to} \;\; C(\R^+, L^1(\R)),$$
\noindent
$$ u^k_0 \rightarrow u^{\infty}_0, \;\; v^k_0 \rightarrow v^{\infty}_0 \;\; \mbox{in}\;\; L^1(\R) \; \mbox{as}\; k \rightarrow \infty,$$
and there exists a continuous function $\omega: \R^+ \to \R^+$ with  $\omega (\mu) \to 0$ as $\mu \to 0$ and 
$$\| \uk_0 (\cdot + \xi) - \uk_0 (\cdot) \|_{L^1(\R)} + \| \vk_0 (\cdot + \xi) - \vk_0 (\cdot) \|_{L^1(\R)} \; \leq \; \omega(|\xi|)\;\; \mbox{for all} \; k>0, \; \xi \in \R.$$

\medskip

\noindent The parameter $k$ is positive  and the interaction function $F: \R^+ \times \R^+ \rightarrow \R^+$ is such that

\medskip

\begin{itemize}
\item[(i)] there exists $\alpha >0$ such that $F \in C^{0, \alpha}(\R^+ \times \R^+)$,
\item[(ii)] $F(u, 0) = F(0, v) =0$ for all $u,v \in \R^+$ and $F(u,v) >0$ for $(u,v) \in (0, \infty) \times (0, \infty)$,
\item[(iii)] $F(\cdot, v)$ and $F(u, \cdot)$ are non-decreasing for all $u,v \in \R^+$.
\end{itemize}

\medskip
\noindent Two cases for $(P_1^k)$ are considered,  when the diffusion coefficients $d_u$
and $d_v$ are both strictly positive (two mobile reactants), and when $d_u>0$
and $d_v=0$ (one mobile and one immobile reactant). 

\bigskip
\noindent The second pair of problems is defined on the half-strip $\st = \{(x,t): 0 < x < \infty, 0 < t < T\}$, and we consider the system
$$
(P_2^k) \left\{\begin{array}{rlrl}\displaystyle
u_t&= d_u u_{xx} -kF(u,v)\quad &\mbox{in}~\st,
\\[.1 in]
\displaystyle v_t&= d_v v_{xx} -kF(u,v)\quad &\mbox{in}~\st,  \\[.1 in]
u(0, t) & = U_0, \quad d_v v_x(0, t) = 0 &\qquad\mbox{for} \quad t \in (0,T), \\[.1 in]
 u(x,0)&=u^k_0(x),\quad v(x,0)=v^k_0(x)   &\qquad\mbox{for} \quad x\in \R,
 \end{array}\right.$$
 \medskip
 
 \noindent where  $u^k_0, v^k_0 \in C^2(\R^+)$ are such that $0 \leq u^k_0 \leq M$, $0 \leq v^k_0 \leq M$ for some $M \geq \max\{ U_0, V_0 \}$, and now 
 $$u^k_0 (x) \rightarrow  0 \;\; \mbox{and} \;\; v^k_0 (x) \rightarrow  V_0 \;\; \mbox{as} \;\; x \rightarrow \infty,$$
 $$\| u^k_0 - u_0^{\infty} \|_{L^1(\R^+)} < \infty,  \;\; \| v^k_0 - v_0^{\infty} \|_{L^1(\R^+)} < \infty,$$
 %$$k \mapsto u_0^k - u_0^{\infty}, \;\; k \mapsto v_0^k - v_0^{\infty} \;\; \mbox{belong to} \;\; C(\R^+, L^1(\R^+)),$$
\noindent 
$$ u^k_0 \rightarrow u^{\infty}_0, \;\; v^k_0 \rightarrow v^{\infty}_0 \;\; \mbox{in}\;\; L^1(\R^+) \; \mbox{as}\; k \rightarrow \infty, $$
 for each $r>0$, there exists a continuous function $\omega_r: \R^+ \to \R^+$ with  $\omega_r (\mu) \to 0$ as $\mu \to 0$ and 
$$\| \uk_0 (\cdot + \xi) - \uk_0 (\cdot) \|_{L^1((r, \infty))} + \| \vk_0 (\cdot + \xi) - \vk_0 (\cdot) \|_{L^1((r, \infty))} \; \leq \; \omega_r(|\xi|)\;\; \mbox{for all} \; k>0, \; |\xi| < r/4,$$

\medskip
\noindent and $k$ and $F$ are as in problem $(P_1^k)$. We  again consider both the case of two mobile reactants, where the diffusion coeffcients $d_u$
and $d_v$ are both strictly positive, and the case of one mobile and one immobile reactant, when $d_u>0$
and $d_v=0$. 

\bigskip
\noindent For each of these four problems, we prove the convergence of solutions $(u^k, v^k)$ on bounded time intervals $(0,T)$ as $ k \rightarrow \infty$ to  a self-similar profile $(u,v)$ in which ${u}$ and ${v}$ are segregated, separated by a free boundary. In each case, the limits ${u}$ of $\uk$ and ${v}$ of $\vk$ are given by the positive and negative parts respectively of a function $w$, that is,
$$u=w^+ \;\;\; \mbox{and} \;\;\; v=-w^-,$$
\noindent where $s^+ = \max\{0, s\}$ and $s^- = \min \{0, s\}$. This limit function $w$ has one of four self-similar forms,   depending on whether $(\uk, \vk)$ satisfy $(P_1^k)$ or $(P_2^k)$, and on whether $d_v>0$ or $d_v=0$. 
 If $(\uk, \vk)$ satisfies $(P_1^k)$, 
there exists a function $f: \R \to \R$ and a constant $a \in \R$ such that
$
w(x,t) = f\left(x/\sqrt{t}\right)$ for $(x,t) \in Q_T$; 
if $d_v>0$, then $a \in \R$ is the unique root of the equation
$  d_u U_0 \int_a^{\infty} e^{\frac{a^2 - s^2}{4d_v}} ds = d_v V_0 \int_{- \infty}^{a} e^{\frac{a^2 - s^2}{4d_u}} ds$,
\noindent and
$$f(\eta) = \left\{ \begin{array}{ll} U_0 \left(  1 - \frac{\int_{- \infty}^{\eta} e^{- \frac{s^2}{4d_u}}\;ds}{\int_{- \infty}^{a} e^{- \frac{s^2}{4d_u}} \;ds}        \right), & \mbox{if} \;\; \eta \leq a,\\
-V_0 \left(  1 - \frac{\int_{\eta}^{\infty} e^{- \frac{s^2}{4d_v}}\;ds}{\int_{a}^{\infty} e^{- \frac{s^2}{4d_v}} \;ds}        \right), & \mbox{if} \;\; \eta>a,
 \end{array} \right.
 $$
\noindent whereas if $d_v=0$, then $a>0$ is the unique root of the equation
$U_0 = \frac{V_0 a}{2d_u} \int_{- \infty}^a e^{\frac{a^2 - s^2}{4 d_u}}\;ds,$
\noindent and
$$
f(\eta) = \left\{ \begin{array}{ll} U_0 \left(  1 - \frac{\int_{- \infty}^{\eta} e^{- \frac{s^2}{4d_u}}\;ds}{\int_{- \infty}^{a} e^{- \frac{s^2}{4d_u}} \;ds}        \right), & \mbox{if} \;\; \eta \leq a,\\
-V_0 , & \mbox{if} \;\; \eta>a.
 \end{array} \right.
$$

\noindent On the other hand, if $(\uk, \vk)$ satisfies $(P_2^k)$, there exists a function $f: \R^+ \to \R$ and a positive constant $a>0$ such that
$w(x,t) = f\left(x/\sqrt{t}\right)$ for $(x,t) \in S_T$; if $d_v>0$, then $a >0 $ is the unique root of the equation
$  d_u U_0 \int_a^{\infty} e^{\frac{a^2 - s^2}{4d_v}} ds = d_v V_0 \int_{0}^{a} e^{\frac{a^2 - s^2}{4d_u}},$
\noindent and
$$
f(\eta) = \left\{ \begin{array}{ll} U_0 \left(  1 - \frac{\int_{0}^{\eta} e^{- \frac{s^2}{4d_u}}\;ds}{\int_{0}^{a} e^{- \frac{s^2}{4d_u}} \;ds}        \right), & \mbox{if} \;\; \eta \leq a,\\
-V_0 \left(  1 - \frac{\int_{\eta}^{\infty} e^{- \frac{s^2}{4d_v}}\;ds}{\int_{a}^{\infty} e^{- \frac{s^2}{4d_v}} \;ds}        \right), & \mbox{if} \;\; \eta>a,
 \end{array} \right.
$$
\noindent whereas if $d_v=0$, then $a>0$ is the unique root of the equation
$U_0 = \frac{V_0 a}{2d_u} \int_{0}^a e^{\frac{a^2 - s^2}{4 d_u}}\;ds,$
\noindent and
$$
f(\eta) = \left\{ \begin{array}{ll} U_0 \left(  1 - \frac{\int_{0}^{\eta} e^{- \frac{s^2}{4d_u}}\;ds}{\int_{0}^{a} e^{- \frac{s^2}{4d_u}} \;ds}        \right), & \mbox{if} \;\; \eta \leq a,\\
-V_0 , & \mbox{if} \;\; \eta>a.
 \end{array} \right.
$$
\medskip

\noindent Clearly, in all four cases, a free boundary is given by the set where $f$ equals zero, which has the form $x=a \sqrt{t}$ where the constant $a$ is determined by a different equation for each problem. Note that only  when $(\uk, \vk)$ satisfies $(P_1^k)$ with $d_u>0$ and $d_v>0$ is the constant $a$ in the corresponding limit problem not necessarily strictly positive, and hence only for this problem is it possible for $v$ to invade $u$ instead of vice versa. Sufficient conditions ensuring $a>0$, $a<0$ or $a=0$ in this case are given in Proposition \ref{propsuffpos}.

\bigskip
\noindent In the last section of the paper, we fix $k$ and initial conditions $u_0$ and $v_0$ such that
$$\| u_0 - u^{\infty}_0 \|_{L^1} < \infty, \;\; \| v_0 - v^{\infty}_0 \|_{L^1} < \infty,$$
\noindent and either
$$u_0(x) \rightarrow U_0, 0 \;\; \mbox{as} \;\;x \rightarrow -\infty, \infty \;\; \mbox{and} \;\; 
v_0(x) \rightarrow 0, V_0 \;\; \mbox{as} \;\;x \rightarrow -\infty, \infty,$$

\noindent in the case of the two full-line problems $(P_1^k)$, or
$$u_0(x) \rightarrow  0 \;\; \mbox{as} \;\;x \rightarrow  \infty \;\; \mbox{and} \;\; 
v_0(x) \rightarrow V_0 \;\; \mbox{as} \;\;x \rightarrow \infty,$$

\noindent in the case of the two half-line problems $(P_2^k)$, 
and then show, by exploiting the $k\to \infty$ results already established, that as $ t \rightarrow \infty$ along a subsequence, $u(\cdot,t)$ and $v(\cdot,t)$ converge, in a certain average sense, to the appropriate one of the same four
self-similar profiles.  The proof uses a scaling argument originally due to Kamin \cite{Kamin}.

\bigskip
\noindent This paper extends the earlier work of \cite{hhp1, hhp2} both  by treating the case of two mobile reactants ($d_u>0$, $d_v>0$) in addition to that of one mobile
reactant ($d_u>0$, $d_v=0$), and in considering the whole-line problem $(P_1^k)$ in addition to the half-line problem $(P_2^k)$. 
 Importantly, we also allow significantly more general initial data than previous work. In \cite{hhp1, hhp2}, the initial
conditions for $(P_2^k)$ are taken to be constant on the half-line $\R^+$, in fact equal to the initial data for the limiting self-similar solution, $u_0^{\infty}|_{\R^+}$, 
$v_0^{\infty}|_{\R^+}$. This implies monotonicity properties in space and time of
solutions $(\uk, \vk)$ of $(P_2^k)$ that are exploited in \cite{hhp1, hhp2} to obtain some compactness  of sequences $\{(\uk, \vk)\}_{k>0}$. 
Here, on the other hand, the initial data $(\uk_0, \vk_0)$ is only supposed to satisfy the hypotheses listed above, and $\uk_0$, $\vk_0$ 
may be non-monotonic in space and can even exceed $U_0$, $V_0$ on parts of the domain. For such initial conditions, monotonicity properties
of $(\uk, \vk)$ are no longer expected, of course, and alternative methods are needed. We exploit some ideas used previously in \cite{CDHMN}, \cite{hkr}
and \cite{hmm}, keeping in mind that here, in contrast to \cite{CDHMN} and \cite{hmm}, our domains are unbounded. Note further that, motivated by the desire to include reaction dynamics of the form $F(u,v) = u^mv^n$ with $0 \leq m <1$, $0 \leq n <1$ (see \cite{Trev2}, for example), we do {not} assume
that $F$ is Lipschitz continuous. Instead, as in \cite{hhp2},  $F$ is assumed to satisfy monotonicity hypotheses   that suffice to establish comparison theorems (see Lemmas \ref{lemcomparisonmain} and \ref{lemcomparisonhalf}) in the absence of  Lipschitz continuity. These monotonicity properties of $F$ also enable the proof of $L^1$-contraction properties (see Lemma \ref{leml1contraction} and Lemma \ref{lemspacetranslatehalf}) giving bounds on differences of space translates, independently of  $d_v$ sufficiently small and of $k$, that yield sufficient compactness to pass to the limits both as $k \to \infty$ and as $d_v \to 0$.

\medskip
\noindent We remark that the form of the self-similar solutions obtained here is clearly due to the presence of the heat operator  and the fact that the same interaction term, $-k F(u,v)$, occurs in each equation in both $(P_1^k)$ and $(P_2^k)$. In  fact, identical limit profiles are obtained for a relatively wide class of interaction terms $-k F(u,v)$ under suitable conditions on $F$,  such as positivity and monotonicity, that suffice to ensure segregation of the two components
and compactness properties of sets of solutions $\{ (\uk, \vk)\}_{k >0}$. 
Interesting potential extensions of this work include investigating possible convergence to other 
types of self-similar solutions  when the diffusion terms $u_{xx}, v_{xx}$ are replaced by nonlinear diffusion terms, and also problems on multi-dimensional spatial domains.

\medskip
\noindent The rest of the paper is organised as follows. In Section 2, we study the whole-line problem $(P_1^k)$, starting with existence and uniqueness of solutions for $(P_1^k)$, first when $d_u >0 $ and $d_v >0$, and then, via some a priori estimates that are also useful in passing to the limit as $k \to \infty$,  when $d_u>0$ and $d_v=0$. A key  bound on $kF(\uk, \vk)$ in $L^1(Q_T)$, independent of $k$ and $d_v \geq 0$, is given in Theorem \ref{lemfbound}. The last part of Section 2 is concerned with the limit of solutions $(\uk, \vk)$ of $(P_1^k)$ as $k \to \infty$, which is characterised as a self-similar solution in Theorem \ref{thmselfsimilar}. This self-similar solution has one of two forms, depending on whether $d_v>0$ or $d_v=0$. Section 3 is devoted to corresponding results for the half-line problem $(P_2^k)$, for which some different arguments are required on account of the boundary at $x=0$. Theorem \ref{lemfboundhalfline} is the half-line counterpart of Theorem \ref{lemfbound}. The two limiting self-similar solutions in this case, one for $d_v>0$ and the other for $d_v=0$,  are given in Theorem \ref{thmselfsimilarhalf}. Finally, in Section 4, the results of the previous sections are used to deduce long-time convergence of solutions of $(P_1^k)$ and $(P_2^k)$ to the appropriate one of the four self-similar  solutions. 

\medskip
\noindent Note that since we are interested in taking limits as $d_v \to 0$, when we write that a given bound is independent of $d_v$, we always mean that the
bound is independent of $d_v \leq D$ for some $D>0$, {\em i.e.} that the bound is independent of $d_v$ sufficiently small. Note also that throughout the paper, our notion of solution of $(P_1^k)$ and $(P_2^k)$ depends on whether $d_v>0$ or $d_v=0$, being classical and weak respectively,
and is made precise in Theorems \ref{lemwholelineexistence}, \ref{thmexistzero}, \ref{lemhalflineexistence} and \ref{thmexistzerohalf} below. 
Various results, such as the comparison principles Lemma \ref{lemcomparisonmain}, \ref{lemcomparisonhalf}, a priori bounds Lemma \ref{lemfbound}, \ref{lemfboundhalfline}, etc.,  hold both when $d_v>0$ and $d_v=0$, with almost identical proofs, and so to avoid duplication, we will present results for $d_v \geq 0$  and understand an appropriate notion of solution in each case. Additionally, we adopt the notational convention that terms multiplied by $d_v$, such as $d_v v_{xx}$, for example, are understood to be simply absent when $d_v=0$.

\medskip

\section{The whole-line case: problem $(P_1^k)$}

\subsection{Existence and uniqueness of solutions for $(P_1^k)$ when $d_u >0 $ and $d_v >0 $}
\label{existuniquewhole}
\noindent We consider first an approximate problem $(P_1^{R, \mu})$ to $(P_1^k)$. Choose $M \geq \max\{ U_0, V_0 \}$, let $R>1$, and consider the problem

$$
(P_1^{R, \mu})\left\{\begin{array}{rlrl}\displaystyle
u_t&= d_u u_{xx} -kF_{\mu}(u,v)\quad &\mbox{in}~(-R, R) \times (0, T),
\\[.1 in]
\displaystyle v_t&= d_v v_{xx} -kF_{\mu} (u,v)\quad &\mbox{in}~(-R, R) \times (0,T), \\[.1 in]
u_x(-R, t) &= u_x(R,t)=0 \quad &\mbox{for} \qquad t \in (0, T),
\\[.1 in]
v_x(-R, t) &= v_x(R,t)=0 \quad &\mbox{for} \qquad t \in (0, T),
\\[.1 in]
 u(x,0)&=u^k_{0,R}(x),\quad v(x,0)=v^k_{0, R}(x)   &\qquad\mbox{for} \quad x\in (-R, R),
\end{array}\right.$$

\noindent where $u^k_{0,R}, v^k_{0, R} \in C^2(\R)$ are such that $0 \leq u^k_{0,R} \leq M$, 
$0 \leq v^k_{0, R} \leq M$  and
\begin{eqnarray}
& & u^k_{0,R} (x) =0  \; \mbox{for}   \; x > \left(1 - {\frac{1}{R}} \right) R, \;\;\; 
  u^k_{0,R} (x) =U_0  \;  \mbox{for} \;  x < - \left(1 - {\frac{1}{R}} \right)R, \label{urini}\\
& & v^k_{0,R} (x) =V_0  \; \mbox{for}   \; x > \left(1 - {\frac{1}{R}} \right) R, \;\;\; 
  v^k_{0,R} (x) =0  \;  \mbox{for} \;  x < - \left(1 - {\frac{1}{R}} \right)R, \label{vrini}
 \end{eqnarray}
 
 \noindent  which defines the functions $u^k_{0,R}, v^k_{0,R}$ on the whole real line. We suppose also that the diffusion coefficients $d_u$ and $d_v$ are both strictly positive. The function $F_{\mu}$ is a regularisation of $F$, such that $F_{\mu} : \R^+ \times \R^+ \rightarrow \R^+$ satisfies
 
 \begin{itemize}
\item[(i)]  $F_{\mu} \in C^{1}(\R^+ \times \R^+)$,
\item[(ii)] $F_{\mu}(u, 0) = F_{\mu}(0, v) =0$ for all $u,v \in \R^+$, and $F_{\mu}(u,v) >0$ for $(u,v) \in (0, \infty) \times (0, \infty)$,
\item[(iii)] $F_{\mu} (\cdot, v)$ and $F_{\mu} (u, \cdot)$ are non-decreasing for all $u,v \in \R^+$,
\item[(iv)] $F_{\mu} \rightarrow F$ in $L^{\infty}_{loc}(\R \times \R)$ as $\mu \rightarrow 0$.
\end{itemize}

\noindent By a solution of $(P_1^{R, \mu})$ we mean a pair $(u,v)$
such that  $u,v \in C^{2,1}(([-R, R] \times [\delta, T] ) \cap C^0([-R, R] \times [0, T])$ for each $\delta >0$
and satisfy $(P_1^{R, \mu})$. 

\begin{lemma}
\label{lemcomparison}
Let $u,v$ and $\tilde{u}, \tilde{v}$ be two solutions of $(P_1^{R, \mu})$ whose initial data satisfy
\begin{equation}
\label{comp1}
u (\cdot, 0) \leq \tilde{u}(\cdot, 0), \;\;\; v (\cdot, 0)  \geq \tilde{v} (\cdot, 0) \;\;\; \mbox{in}\;\; (-R, R).
\end{equation}

\noindent Then 
\begin{equation}
\label{comp2}
u(\cdot, t) \leq \tilde{u}(\cdot, t), \;\;\; v(\cdot, t) \geq \tilde{v}(\cdot, t) \;\;\; \mbox{in} \;\; (-R, R) \times (0,T).
\end{equation}
\end{lemma}

\medskip
\pf This follows from \cite[p 241, Lem. 5.2 and p 244, Thm. 5.5]{Vol3}   applied to the new system obtained from
$(P_1^{R, \mu})$ under the change of variables $u \mapsto u$ and $v \mapsto V_0 - v$ (note that in the notation of \cite{Vol3}, $u = (u_1, u_2)$ is a vector). \qed

\bigskip
\noindent The following corollary is immediate from Lemma \ref{lemcomparison}. 

\begin{corollary}
\label{corunique}
For given initial data $u^k_{0,R}$, $v^k_{0,R}$, there is at most one solution $(u^k_{R, \mu},v^k_{R, \mu})$ of  $(P_1^{R, \mu})$. 
\end{corollary}

\medskip
\noindent We also have the following bound, which is easily proved using the scalar maximum principle. 

\begin{lemma}
\label{lembound}
Let $(u^k_{R, \mu},v^k_{R, \mu})$ be a solution of $(P_1^{R, \mu})$. Then
\begin{equation}
\label{bound}
0 \leq u^k_{R, \mu} \leq M, \;\; 0 \leq v^k_{R, \mu} \leq M \;\;\; \mbox{on}\;\;\; (-R, R) \times (0,T). 
\end{equation}
\end{lemma}

\pf  We define
\begin{eqnarray*}
\mathcal{L}_1 (u) & := & u_t - d_u  u_{xx}  + kF(u,v), \\
\mathcal{L}_2 (v) & := & v_t - d_v  v_{xx} +  kF(u,v).
\end{eqnarray*}
Since $\mathcal{L}_i (0) =0$ and $\mathcal{L}_i(M) \geq 0$ for
$i=1,2$, the assertion follows from the maximum  principle. \qed

\bigskip

\begin{lemma}
\label{lemaux1}
There exists a unique solution $(u^k_{R, \mu},v^k_{R, \mu})$ of $(P_1^{R, \mu})$.
 \end{lemma}
 
 \pf It follows from Lunardi \cite[Prop. 7.3.2]{lunardi} that there exist  $u^k_{R, \mu}, v^k_{R, \mu}$ and $T^* \in (0, T]$ such that
 $u^k_{R, \mu}, v^k_{R, \mu} \in C^{2,1}(([\delta, T] \times [-R, R]) \cap C^0([-R, R] \times [0, T])$ for each $\delta >0$
and satisfy $(P_1^{R, \mu})$ with $T$ replaced by $T^*$. That in fact we can take $T^* = T$ is a consequence of Lemma \ref{lembound}, and uniqueness of the solution is given by Corollary \ref{corunique}. \qed

\bigskip

\noindent We now introduce  a class  of cut-off functions. \textcolor{black}{First define an even, non-negative cut-off function $\psi^1 \in C^{\infty}(\R)$ such that  $0 \leq \psi^1(x) \leq 1$ for all $x \in \R$, $\psi^1(x) = 1$ when $|x| \leq 1$,  and $\psi^1(x)=0$ when $|x| \geq 2$. Then given $L \geq 1$, define  the family of cut-off functions $\psi^L \in C^{\infty}(\R)$ by $\psi^L(x) = 1$ when $|x| \leq L$ and $\psi^L(x) = \psi^1(|x| + 1 - L)$ when $|x| \geq L$. Clearly $\psi^L$, $\psi^L_x$ and $\psi^L_{xx}$ are bounded in $L^{\infty}(\R)$ independently of $L$, and $\psi^L_x$ and $\psi^L_{xx}$, being supported on sets of measure at most two, are also bounded in $L^1(\R)$ independently of $L$.} Let $Q_{L, T}$ denote the truncated space-time domain  $(-L, L) \times (-T, T)$. In the following, 
$C(L)$ denotes some $L$-dependent constant  which varies according to context.

\begin{lemma}
\label{lemkuvest}
Let  $L>0$. Then there exists a constant $C(L)$ such that if $R>L+1$, then
\begin{equation}
\label{kuvest}
k \iint_{Q_{L, T}} F_{\mu}(u^k_{R, \mu}, v^k_{R, \mu})  ~ dxdt \;  \leq \; C(L),
\end{equation}

\noindent for all $k, \mu >0$.
\end{lemma}

\pf Multiplying the equation for $u^k_{R, \mu}$ in $(P_1^{R, \mu})$ by $\psi^L$ and integrating over $Q_{L+1, T}$
gives that
$$\int_{-L-1}^{L+1} \psi^L \{ u^k_{R, \mu}(\cdot, T) - u^k_{0,R}(\cdot) \} \; dx = d_u \iint_{Q_{L+1, T}} u^k_{R, \mu}
\psi^L_{xx}\; dx dt - k \iint_{Q_{L+1, T}} F_{\mu}(u^k_{R, \mu}, v^k_{R, \mu}) \psi^L ~ dxdt,$$

\noindent which, together with Lemma \ref{lembound} and the definition of $\psi^L$,  yields (\ref{kuvest}). \qed 

\begin{lemma}
\label{lemuvhone}
The solutions $u^k_{R, \mu}, v^k_{R, \mu}$ are bounded  in $L^2(0, T; H^1_{loc}(\R))$ independently of
$k, R, \mu$.
\end{lemma}

\pf 
We prove the bound for $u^k_{R, \mu}$. Suppose that $R>L+1$. Then multiplying the equation for $u^k_{R, \mu}$ by $u^k_{R, \mu} \psi^L$ and integrating over $Q_{L+1, T}$ gives that
\begin{eqnarray*}
\lefteqn{\frac{1}{2} \int_{-L-1}^{L+1} \psi^L \; \{ u^k_{R, \mu}(\cdot, T)^2 - u^k_{0,R}(\cdot)^2 \} \; dx \leq 
- d_u \iint_{Q_{L+1, T}} (u^k_{R, \mu})_x^2 \; \psi^L \; dx dt \;} \\ & & \hspace*{9.5cm} +\frac{d_u}{2}  \iint_{Q_{L+1, T}} (u^k_{R, \mu})^2
\psi^L_{xx} \; dx dt,
\end{eqnarray*}
\noindent since $F_{\mu}(u^k_{R, \mu}, v^k_{R, \mu}) \geq 0$. The result again follows using
Lemma \ref{lembound} and the definition of $\psi^L$. \qed

\medskip
\noindent In order to prove that the sets $\{ u^k_{R, \mu} : k, R, \mu >0 \}$, $\{ v^k_{R, \mu} : k, R, \mu >0 \}$ are each relatively compact in $L^2_{loc}(\R \times (0,T))$, we first give estimates of the differences of space and time translates of $u^k_{R, \mu}$ and  $v^k_{R, \mu}$. 

\begin{lemma}
\label{lemspace}
For each $L>0$, there exists a constant $C(L)$ such that
\begin{eqnarray*}
& & \iint_{Q_{L, T}} ( u^k_{R, \mu}(x + \xi, t) - u^k_{R, \mu}(x,t))^2\; dx dt \leq C(L) | \xi|^2,\\
& & \iint_{Q_{L, T}} ( v^k_{R, \mu}(x + \xi, t) - v^k_{R, \mu}(x,t))^2\; dx dt \leq C(L) | \xi|^2,\\
\end{eqnarray*}

\noindent for all $\xi \in \R$, $|\xi| \leq L$. 

\end{lemma}

\pf As a result of the gradient bounds in Lemma \ref{lemuvhone}, the proof of this closely follows the proof of \cite[Lemma 2.6]{CDHMN} and we omit the details. \qed

\begin{lemma}
\label{lemtime}
For each $L>0$, there exists a constant $C(L)$ such that
\begin{eqnarray*}
& & \iint_{Q_{L, T- \tau}} ( u^k_{R, \mu}(x , t + \tau) - u^k_{R, \mu}(x,t))^2\; dx dt \leq C(L) \tau,\\
& & \iint_{Q_{L, T- \tau }} ( v^k_{R, \mu}(x, t + \tau) - v^k_{R, \mu}(x,t))^2\; dx dt \leq C(L) \tau,\\
\end{eqnarray*}

\noindent for all $\tau \in (0,T)$. 

\end{lemma}

\pf The gradient bounds in Lemma \ref{lemuvhone} together with Lemma \ref{lemkuvest} enable the proof of \cite[Lemma 2.7]{CDHMN} to be easily adapted. \qed

\bigskip
\noindent We can now establish the existence of a classical solution of the original problem $(P_1^k)$ on $Q_T$
when both diffusion coefficients $d_u$ and $d_v$ are strictly positive. 

\begin{theorem}
\label{lemwholelineexistence}
Suppose that $d_u>0$ and $d_v>0$. Then given $k>0$, there exists a classical solution $(\uk, \vk)$ of $(P_1^k)$ such that
for each $\delta>0$, $J>0$ and $p \geq 1$,
\begin{equation}
\label{reg-pone}
\uk, \vk \in C^{2,1}(\R \times [\delta, T] ) \cap C^0(\R \times [0, T]) \cap W^{2,1}_p((-J, J) \times (0, T)),
\end{equation}
and
\begin{equation}
\label{linfbound}
0 \leq u^k \leq M, \;\;\; 0 \leq v^k \leq M \;\;\; \mbox{on} \;\;\; \R \times (0, T).
\end{equation}
\end{theorem}

\pf  Let  $u^k_{0,R}, v^k_{0, R}$ be
as in the formulation of problem $(P_1^{R, \mu})$ and such that as $ R \rightarrow \infty$, 
$u^k_{0,R} \rightarrow u^k_0$ and $v^k_{0,R} \rightarrow v^k_0$ in $C^1_{loc}(\R)$. 
Then given
$R_n \rightarrow \infty$ and $\mu_n \downarrow 0$, it follows from the Fr\'echet-Kolmogorov Theorem (see, for example, \cite[Corollary 4.27]{brezisnew}) and Lemmas \ref{lembound}, \ref{lemspace} and \ref{lemtime},
that there exist subsequences 
 $\{R_{n_j} \}_{j=1}^{\infty}$, $\{{\mu}_{n_j} \}_{j=1}^{\infty}$
 and functions $u^k \in L^{\infty}(Q_T)$ and $v^k \in L^{\infty}(Q_T)$ such that
$$u^{k}_{R_{n_j}, \mu_{n_j}} \rightarrow u^k, \;\;\; v^{k}_{R_{n_j}, \mu_{n_j}} \rightarrow v^k \;\;\; 
\mbox{strongly in}\; 
L^2_{loc}(Q_T) \;\; \mbox{and} \;\; a.e. \; \mbox{in} \; Q_T,$$

\noindent as $j \rightarrow \infty$. We can then easily pass to the limit in the weak form of $(P_1^k)$.
To see that the solution is  in fact classical, note first that for a fixed $k$, the term $kF(u^k,v^k)$
is in $L^{\infty}(Q_T)$, which, since $u_0^k, v_0^k \in C^2(\R)$,  implies that $\uk, \vk \in W^{2,1}_p((-J, J) \times (0, T))$ for each $J>0$ and $p \geq 1$, and hence $u^k, v^k \in C^{1 + \lambda, \frac{1 + \lambda}{2}} (\R \times [0, T])$ for each $\lambda \in (0,1)$. Since $F \in C^{0, \alpha}(\R^+ \times \R^+)$, it then follows that $F(u^k, v^k)$ is H\"older continuous and so $u^k, v^k \in C^{2 + \lambda, \frac{2 + \lambda}{2}}(\R \times (0, T])$ for some $\lambda > 0$. The bounds (\ref{linfbound}) are immediate from Lemma \ref{lembound}. \qed

\bigskip \medskip

\noindent To show uniqueness, we use the following comparison theorem for $(P_1^k)$, proved with arguments inspired by \cite[Lemma 2.7]{hhp2}. Note that this result covers  both the case $d_u >0$, $d_v>0$ and the case $d_u>0$, $d_v =0$, and  the monotonicity properties of $F$ are exploited  to overcome the fact that $F$ is not assumed to be Lipschitz continuous. For an alternative approach when $F$ is Lipschitz and $d_v=0$, see \cite[Lemma 5]{hmm}. 

\begin{lemma}
\label{lemcomparisonmain}
Suppose that $d_u > 0$, $d_v \geq 0$, and let $(\ou, \ov)$, $(\uu, \uv)$ be such that  for each $J>0$ and $p \geq 1$, $\ou, \uu\in L^{\infty}(Q_T) \cap  W^{2,1}_p((-J, J) \times (0, T))$, $\ov, \uv \in L^{\infty}(Q_T)  \cap W^{2,1}_p((-J, J) \times (0, T))$ if $d_v >0$,   $\ov, \uv \in  L^{\infty}(Q_T) \cap W^{1, \infty}(0, T; L^{\infty}((-J, J)))$ if $d_v=0$, and $(\ou, \ov)$, $(\uu, \uv)$ satisfy
$$\begin{array}{l} \ou_t \geq d_u \ou_{xx} - kF(\ou, \ov),\\ \ov_t \leq d_v \ov_{xx} - kF(\ou, \ov),\end{array} \qquad 
\begin{array}{l} \uu_t \leq d_u \uu_{xx} - kF(\uu, \uv),\\ \uv_t \geq d_v \uv_{xx} - kF(\uu, \uv),\end{array} \qquad \mbox{in} \;\; Q_T, 
$$
and
$$\ou(\cdot, 0) \geq \uu(\cdot,0), \qquad \ov(\cdot,0) \leq \uv(\cdot, 0) \quad  \mbox{on}\;  \; \R.$$

\medskip
\noindent Then 

$$
\ou \geq \uu \quad \mbox{and} \quad \ov \leq \uv \quad \mbox{in} \;\; Q_T.
$$
\end{lemma}
\medskip

\pf Let $u:= \uu-\ou$, $v:= \ov - \uv$, $u_0:= \uu(\cdot, 0) - \ou(\cdot, 0)$ and $v_0 := \ov(\cdot, 0) -
\uv(\cdot, 0)$. Then
\begin{eqnarray}
u_t & \leq & d_u u_{xx} - k \{ F(\uu, \uv) - F(\ou, \ov) \} \quad \mbox{ in} \;\; Q_T, \label{compa1}\\
v_t  & \leq & d_v v_{xx} -  k \{ F(\ou, \ov) - F(\uu, \uv) \} \quad \, \mbox{ in} \;\; Q_T, \label{compa2}
\end{eqnarray}

\noindent and
$$u_0 \leq 0, \quad v_0 \leq 0 \;\; \mbox{on} \;\; \R.$$

\noindent Now take  a smooth non-decreasing convex function $m^+: \R \to \R$ with
$$m^+ \geq 0, \quad m^+(0) =0, \quad (m^+)'(0)=0, \quad m^+(r) \equiv 0 \;\;\; \mbox{for} \;\; r \leq 0, \quad m^+(r) = |r| - \frac{1}{2} \quad \mbox{for} \;\; r >1,$$

\noindent and for each $\alpha >0$, define the functions
$$\ma^+ (r) : = \alpha m^+ \left( \frac{r}{\alpha} \right),$$

\noindent which  as $\alpha \to 0$ approximate the positive part of $r$. Then multiplying (\ref{compa1}) by $(\ma^+)'(u)$ and (\ref{compa2}) by  $(\ma^+)'(v)$ gives
\begin{eqnarray*}
(\ma^+)'(u) u_t & \leq & d_u (\ma^+)'(u) u_{xx} - k (\ma^+)'(u) \{ F(\uu, \uv) - F(\ou, \ov) \} \quad \mbox{in} \;\; Q_T,\\
(\ma^+)'(v) v_t & \leq & d_v (\ma^+)'(v) v_{xx} - k (\ma^+)'(v) \{ F(\ou, \ov) - F(\uu, \uv) \} \quad \, \mbox{in} \;\; Q_T,
\end{eqnarray*}

\noindent and it follows from adding these inequalities that
\begin{eqnarray}
\lefteqn{ (\ma^+)'(u) u_t + (\ma^+)'(v) v_t \leq d_u (\ma^+)'(u) u_{xx} + d_v (\ma^+)'(v) v_{xx}} \label{compa3}
\\
& & \hspace*{6cm} -k [ (\ma^+)'(v) - (\ma^+)'(u)]  \{ F(\ou, \ov) - F(\uu, \uv) \}. \nonumber
\end{eqnarray}

\noindent Now with  $\psi^L$  the cut-off functions defined before Lemma \ref{lemkuvest}, multiplying by $\psi^L$ and integrating over $\R \times (0, t_0)$, $t_0 \in (0, T]$, gives  
\begin{eqnarray*}
 \int_0^{t_0}\int_{\R} u_{xx} (\ma^+)'(u) \psi^L~dxdt  & =  & -  \int_0^{t_0}\int_{\R} u_x [ (\ma^+)''(u) u_x \psi^L + (\ma^+)'(u) \psi^L_x]~dxdt \\
&  \leq  & -  \int_0^{t_0}\int_{\R} u_x (\ma^+)'(u) \psi^L_x~dxdt =  \int_0^{t_0}\int_{\R} \ma^+(u) \psi^L_{xx}~dxdt,
\end{eqnarray*}

\noindent \textcolor{black}{since $(\ma^+)''(u) \geq 0$ because $\ma^+$ is convex.} So  (\ref{compa3})  yields
 \begin{eqnarray*}
\lefteqn{\int_{\R} \psi^L(x) [\ma^+(u) + \ma^+(v)](x, t_0)~dx  \leq  \int_{\R} \psi^L(x) [\ma^+(u) + \ma^+(v)](x, 0)~dx} \\ & & \hspace*{4cm}   + \int_0^{t_0}  \int_{\R} \psi^L_{xx}(x)\{ d_u \ma^+(u) + d_v \ma^+(v) \} ~dxdt  \\ &  & \hspace*{4cm} - k \int_0^{t_0} \int_{\R} \psi^L(x) [ (\ma^+)'(v) - (\ma^+)'(u) ]  \{ F(\ou, \ov) - F(\uu, \uv)\}~dxdt, 
\end{eqnarray*}

\noindent and  letting $\alpha \to 0$ gives
 \begin{eqnarray*}
\lefteqn{\int_{\R} \psi^L(x) [u^+ + v^+](x, t_0)~dx  \leq  \int_{\R} \psi^L(x) [u^+ + v^+](x, 0)~dx  + \int_0^{t_0}  \int_{\R} \psi^L_{xx}(x)\{ d_u u^+ + d_v v^+ \} ~dxdt}  \\ &  & \hspace*{4cm} - k \int_0^{t_0} \int_{\R} \psi^L(x) [ (\mbox{sgn} v)^+ - (\mbox{sgn} u)^+ ]  \{ F(\ou, \ov) - F(\uu, \uv)\}~dxdt, 
\end{eqnarray*}

\noindent where $u^+:= \max(u, 0)$. Then $(u^+ + v^+)(\cdot, 0) = 0$, and the expression
$$\sharp : = [ (\mbox{sgn} v)^+ - (\mbox{sgn} u)^+ ]  \{ F(\ou, \ov) - F(\uu, \uv)\}$$ \noindent is non-zero only if either
\begin{itemize}
\item[(i)] $(\mbox{sgn} v)^+ =1$ and $(\mbox{sgn} u)^+ =0$, in which case $\ov \geq \uv$ and $\uu \leq \ou$, so that $ F(\ou, \ov) - F(\uu, \uv) \geq 0$, because $F(\cdot, v)$ and $F(u, \cdot)$ are non-decreasing for all $u,v \in \R^+$, and hence $\sharp \geq 0$, or
\item[(ii)] $(\mbox{sgn} v)^+ =0$ and $(\mbox{sgn} u)^+ =1$, in which case $\ov \leq \uv$ and $\uu \geq \ou$, so that $ F(\ou, \ov) - F(\uu, \uv) \leq 0$, and hence, again,  $\sharp \geq 0$.
\end{itemize}

\noindent Thus 
\begin{equation}
\label{complast}
\int_{\R} \psi^L [u^+ + v^+]( x, t_0)~dx \leq  \int_0^{t_0} \int_{\R} [d_u u^+ + d_v v^+] | \psi^L_{xx}| \, dxdt.
\end{equation}

\noindent Now the right-hand side of (\ref{complast}) is bounded independently of $L$. So by Lebesgue's monotone convergence theorem, $u^+, v^+ \in L^{\infty}(0, T; L^1(\R))$, and thus the right-hand side of  (\ref{complast}) in fact tends to $0$ as $L \to \infty$. Hence
$$[u^+ + v^+](\cdot, t_0) =0 \quad \mbox{on} \;\;  \R,$$

\noindent and the result follows. \qed

\bigskip

\noindent The following corollary is immediate from Lemma \ref{lemcomparisonmain}. 

\begin{corollary}
\label{cor1}
Suppose $d_u > 0$ and $d_v > 0$. Then given $k>0$, there is at most one solution $(u^k, v^k)$ of $(P_1^k)$ in $L^{\infty}(Q_T) \cap W^{2,1}_p((-J, J) \times (0, T))$ for all $J>0$, $p \geq 1$. 
\end{corollary}

\bigskip

\subsection{Existence and uniqueness of solutions for $(P_1^k)$ when $d_u>0$ and $d_v=0$}
\label{existuniquenodiffusion}
\noindent Next we prove some preliminary estimates that  will be used in the following both to prove existence of solutions of  $(P_1^k)$ when $d_u>0$ and $d_v=0$, and, in the next section,  to study the limit of  $(P_1^k)$ as $k \to \infty$.

\bigskip 
\noindent The following bound for $kF(\uk, \vk)$ is key. Note that  $kF(\uk, \vk)$ is controlled by $\uk$ on part of the spatial domain and by $\vk$ on the other part, \textcolor{black}{due to the fact that $u_0^k$ is bounded in $\| \cdot \|_{L^1(\R^+)}$ independently of $k$, and  $v^k_0$ is bounded in $\| \cdot \|_{L^1(\R^-)}$  independently of $k$.}  A similar phenomenon occurs in the proof of the corresponding estimate in the half-line case, Lemma \ref{lemfboundhalfline}, in which  $kF(\uk, \vk)$ is controlled by $\uk$ on $(1, \infty) \times (0, T)$,  and by $\vk$ in the boundary region $(0, 1) \times (0, T)$. \smallskip

\begin{lemma}
\label{lemfbound}
There exists a constant $C>0$, independent of  $d_v \geq 0$ and $k>0$,  such that for  any solution $(u^k, v^k)$ of $(P_1^k)$
satisfying (\ref{linfbound}), we have
$$
\int_0^T \int_{\R} k F(u^k, v^k) ~ dx\,dt \leq C.
$$
\end{lemma}

\pf \textcolor{black}{Define a cut-off function $\phi^1 \in C^{\infty}(\R)$ such that $0 \leq \phi^1(x) \leq 1$ for all $x \in \R$, $\phi^1=1$ when $x \in [0,1]$, and $\phi^1(x)=0$ when $x \not\in (-1, 2)$. Then given $L \geq 1$, define $\phi^L \in C^{\infty}(\R)$ by $\phi^L(x) = \phi^1(x)$ if $x \leq 0$,  $\phi^L(x) = 1$ when $x \in [0,L]$, and $\phi^L(x) = \phi^1(x+1-L)$ when $x \geq L$, and define $\tilde{\phi}^L \in C^{\infty}(\R)$ by $\tilde{\phi}^L(x) = \phi^L(-x)$ for all $x \in \R$. Note that $0 \leq \phi^L(x), \tilde{\phi}^L(x)  \leq 1 $ for all $x\in \R$,  
  and $\phi^L_x$,  $\phi^L_{xx}$,  $\tilde{\phi}^L_x$ and $\tilde{\phi}^L_{xx}$ are bounded in both $L^{\infty}(\R)$ and $L^1(\R)$ independently of $L$. } Consider first the case when $d_u>0$ and $d_v >0$. Then multiplying the equation for $u^k$ by $\phi^L$ and integrating over $\R \times (0, t_0)$, $t_0 \in (0, T]$, gives that 
\begin{eqnarray} \nonumber \lefteqn{\int_{-1}^{\infty} \phi^L(x)  u^k(x,t_0)\,dx +  \int_0^{t_0} \int_{-1}^{\infty} \phi^L(x) k F(u^k, v^k) \, dxdt  =} \\ & &  \hspace*{5cm}d_u \int_0^{t_0} \int_{-1}^{\infty} \phi^L_{xx}(x) u^k(x,t)\,dxdt  + \int_{-1}^{\infty} \phi^L(x)  u_0^k(x)\,dx, \label{fboundeq1}
\end{eqnarray}

\noindent which, since the definition of $\phi^L$ and  the facts that $0 \leq u^k \leq M$ and  $\| u_0^k \|_{L^1(\R^+)}$ is bounded independently of $k$ imply that the right-hand side of \eqref{fboundeq1} is bounded independently of $L$ and $k$, gives the existence of $C>0$ such that for all $k >0$ and $t_0 \in (0, T]$, 
\begin{equation}
  \int_{-1}^{\infty} \phi^L(x) u^k(x,t_0)\,dx + \int_0^{t_0} \int_{-1}^{\infty} \phi^L(x) k F(u^k, v^k) \, dxdt  \; \leq \; C,
\label{forl1bound1}
\end{equation}
\noindent and then, since $u^k \geq 0$, letting $L \to \infty$ using Lebesgue's monotone convergence theorem gives 
\begin{equation}
\label{festone}
\int_0^T \int_0^{\infty} kF(u^k, v^k) \, dxdt \leq C. 
\end{equation}

\noindent  (Note that if we  had $d_u=0$ instead of $d_u>0$, then (\ref{forl1bound1})  could be proved
likewise,  with the  first term on the right-hand side of  (\ref{fboundeq1}) absent due to the lack of diffusion term.)  

\medskip
\noindent Similarly, since $\|v_0^k\|_{L^1(\R^-)}$ is bounded independently of $k$,
multiplying the equation for $v^k$ by  $\tilde{\phi}^L$ and  integrating over $\R \times (0, t_0)$ yields that $C$ can be chosen large enough that for all $L, k >0$, we also have
\begin{equation}
  \int_{-\infty}^{1} \tilde{\phi}^L(x) v^k(x,t_0)\,dx + \int_0^{t_0} \int_{-\infty}^{1} \tilde{\phi}^L(x) k F(u^k, v^k) \, dxdt  \; \leq \; C,
\label{forl1bound2}
\end{equation}
and hence, since $v^k \geq 0$, letting $L \to \infty$ yields that 
\begin{equation}
\label{festtwo}
\int_0^T \int_{- \infty}^{0} k F(u^k, v^k)\, dx \leq C.
\end{equation}
 The result then follows from (\ref{festone}) and (\ref{festtwo}). 
\qed

\bigskip
\begin{lemma}
\label{lemloneboundwholeline}
 There exists a constant $C>0$, independent of  $d_v \geq 0$ and $k>0$,  such that for all $k>0$ and  any solution $(u^k, v^k)$ of $(P_1^k)$
satisfying (\ref{linfbound}),
$$
\| u^k(\cdot,t) - u_0^{\infty} \|_{L^1{(\R)}}~ \leq C \qquad \mbox{and} \qquad  \| v^k(\cdot,t) - v_0^{\infty} \|_{L^1{(\R)}}
 \leq C \;\;\; \mbox{for all} \;\;t \in [0, T]. 
$$

\end{lemma}

\pf  Note first that it follows immediately from \eqref{forl1bound1}, \eqref{forl1bound2} and Lebesgue's monotone convergence theorem that there exists $C>0$, independent of  $d_v \geq 0$ and $k>0$,
such that  
\begin{equation}
\label{l1bound1}
\int_0^{\infty} u^k(x,t_0) ~dx \leq C \quad \mbox{and}  \quad \int_{-\infty}^0 v^k(x,t_0) ~ dx \leq C \quad \mbox{for all} \;\; t_0 \in [0, T]. 
\end{equation}

\noindent Now choose a smooth convex function $m: \R \to \R$ with
$$m \geq 0, \qquad m(0) =0, \qquad m'(0)=0, \qquad m(r) = |r| - \frac{1}{2} \qquad \mbox{for} \;\; |r| >1,$$

\noindent and for each $\alpha >0$, define the functions
$$\ma (r) : = \alpha m\left( \frac{r}{\alpha} \right),$$

\noindent which approximate the modulus function as $\alpha \to 0$, and define
$\huk : = U_0 - u^k$. Then 
$$
\begin{array}{rlrl}\displaystyle
\huk_t&= d_u \huk_{xx} + kF(u^k,v^k)\quad &\mbox{in}~\qt
\\[.1 in]
 \huk (x,0)&= U_0  - u^k_0(x), &\qquad\mbox{for} \quad x\in \R.
\end{array}$$

\noindent  Now with $\tpl$ as in Lemma \ref{lemfbound}, 
\begin{eqnarray*}
 \int_{\R} \huk_{xx} \ma'(\huk) \tpl~dx  & =  & - \int_{\R} \huk_x [ \ma''(\huk) \huk_x \tpl + \ma'(\huk) \tpl_x]~dx \\
&  \leq  & - \int_{\R} \huk_x \ma'(\huk) \tpl_x~dx = \int_{\R} \ma(\huk) \tpl_{xx}~dx,
\end{eqnarray*}

\noindent so multiplying the equation for $\huk$ by $\tpl \ma'(\huk)$ and integrating over 
$\R \times (0, t_0)$, $t_0 \in (0, T)$,  gives that
\begin{eqnarray}
 \lefteqn{\;\;\;\int_{\R} \tpl \ma(\huk(x, t_0))~dx \leq \int_{\R} \tpl \ma(\huk(x,0))~dx
\,+ \, d_u  \int_0^{t_0} \int_{\R} \ma(\huk) \tpl_{xx} ~ dxdt} \nonumber
\\ & & \hspace*{8cm}  \,+\, \int_0^{t_0} \int_{\R}
kF(u^k, v^k) \ma'(\huk) \tpl~dxdt, \label{prelim}
\end{eqnarray}

\noindent and then letting $\alpha \to 0$ in (\ref{prelim}) yields 
\begin{eqnarray}
\lefteqn{\int_{\R} \tpl |\huk(x, t_0)| \,dx \leq}\nonumber  \\ & & \hspace*{1cm} \int_{\R} \tpl |\huk(x,0)| \, dx
+  d_u  \int_0^{t_0} \int_{\R} |\huk| \tpl_{xx} \, dxdt
 + \int_0^{t_0} \int_{\R}
kF(u^k, v^k) \mbox{sgn}(\huk) \tpl \, dxdt.\label{alphazero}
\end{eqnarray}

\noindent Now by Lemma \ref{lemfbound}, (\ref{linfbound}), and the fact that 
$\| u_0^k - u_0^{\infty}\|_{L^1(\R)}$ is bounded independently of $k$,  the right-hand side of (\ref{alphazero}) is bounded independently of $L$ and $k$. So it follows from (\ref{alphazero}) that there exists $C$, independent of $k$, such that 
\begin{equation}
\label{ing1} \int_{-\infty}^0 | u^k(x, t_0) - U_0|\, dx  \, \leq \, C \quad \mbox{for all} \;\;t_0 \in (0, T). 
\end{equation}

\noindent Then taking $\phi^L$ as in Lemma \ref{lemfbound},  multiplying the equation satisfied by $\hat{v}^k:= V_0 - v^k$ by $\phi^L \ma(\hat{v}^k)$ and again integrating over $\R \times (0, t_0)$ gives, using a similar argument to above, that $C$ can be chosen large enough that we also have that 
\begin{equation}
\label{ing2}
\int_0^{\infty} | v^k(x, t_0) - V_0|\, dx \, \leq \,C \quad \mbox{for all} \;\;t_0 \in (0, T). 
\end{equation}

\noindent  The result follows from (\ref{ing1}), (\ref{ing2}), and \eqref{l1bound1}. \qed

\bigskip

\noindent We prove next a bound for the $L^2$-norm of the space derivatives $u_x$ and $v_x$. 

\medskip
\begin{lemma}
\label{leml2derivbound}
Suppose that $d_u >0$  and $d_v \geq 0$. Then there exists a constant $C$, independent of  $d_v \geq 0$  and $k>0$, 
such that for   any solution $(u^k, v^k)$ of $(P_1^k)$
satisfying (\ref{linfbound}),
\begin{equation}
\label{derivbounds}
d_u \int_0^T \int_{\R} (u^k_x)^2(x,t) ~dxdt \leq C \quad \mbox{and} \quad  d_v \int_0^T \int_{\R} (v^k_x)^2(x,t)~dxdt \leq C.
\end{equation}
\end{lemma}

\pf Let $\phi^L$ and $\tpl$ be as in the proof of Lemma \ref{lemfbound}. Then multiplication of the equation for $u^k$ by $u^k \phi^L$ and integration over $Q_T$ gives
\begin{eqnarray}
\lefteqn{\frac{1}{2} \int_{\R} \phi^L(x) (u^k)^2 (x, T)~dx + d_u \int_0^T \int_{\R} (u^k_x)^2 \phi^L(x)~dxdt =}\nonumber \\ & & \frac{1}{2} \int_{\R} \phi^L(x) (u^k)^2 (x, 0)~dx  + \frac{d_u}{2} \int_0^T \int_{\R} (u^k)^2(x,t) \phi^L_{xx}(x) ~ dxdt -  \int_0^T\int_{\R} k u^k F(u^k, v^k) \phi^L(x)~dxdt  \nonumber \\
& & \hspace*{2cm} \leq \; \frac{1}{2} \int_{\R} \phi^L(x) (u^k)^2 (x, 0)~dx  + \frac{d_u}{2} \int_0^T \int_{\R} (u^k)^2(x,t) \phi^L_{xx}(x) ~ dxdt, \label{lastline1}
\end{eqnarray}

\noindent since $F(u^k, v^k) \geq 0$. Now it follows from (\ref{linfbound}) and the definition of $\phi^L$ that the right-hand side of (\ref{lastline1}) is bounded independently of $L$. It thus follows that $d_u \int_0^T \int_{\R} (u^k_x)^2 \phi^L~dxdt$ is bounded independently of $L$ and hence, using Lebesgue's monotone convergence theorem to let $L \to \infty$ in $ \int_0^T \int_0^{\infty} (u^k_x)^2 \phi^L~dxdt$, that letting $L \to \infty$ in (\ref{lastline1}) implies that there exists a constant $C_1>0$ such that
$$   d_u \int_0^T \int_0^{\infty} (u^k_x)^2~ ~dxdt \leq C_1 + \frac{1}{2} \int_0^{\infty} (u^k)^2(x,0)~dx \leq \frac{C}{2}, $$
\noindent where the constant $C$ is independent of $d_u, k>0$, by  (\ref{linfbound}) and the fact that $\| u_0^k - u_0^{\infty}\|_{L^1{(\R)}}$ is bounded independently of $k$.  Then letting $\huk : = U_0 - u^k$, multiplying the equation for $\huk$ by $\huk \tpl$ and integrating over $Q_T$ yields
\begin{eqnarray}
\lefteqn{\;\;\; \frac{1}{2} \int_{\R} \tpl (\huk)^2 (x, T)~dx + d_u \int_0^T \int_{\R} (\huk_x)^2 \tpl~dxdt =} \label{lastline2} \\ & & \frac{1}{2} \int_{\R} \tpl (\huk)^2 (x, 0)~dx  + \frac{d_u}{2} \int_0^T \int_{\R} (\huk)^2(x,t) \tpl_{xx} ~ dxdt -  \int_0^T\int_{\R} k \huk F(u^k, v^k) \tpl~dxdt.  \nonumber
\end{eqnarray}

\noindent Now $\huk$ may not be non-negative, but we can call on Lemma \ref{lemfbound} to deduce that the right-hand side of (\ref{lastline2}) is bounded independently of $L$, so that arguing similarly to before gives that $C_1$ can be chosen larger if necessary so that
$$ d_u \int_0^T \int_{- \infty}^0 (u^k_x)^2~ ~dxdt \; \leq \; C_1 + \frac{1}{2}\int_{- \infty}^1 (\huk)^2(x,0)~dx
\; + \; M  \int_0^T\int_{\R} k  F(u^k, v^k) ~dxdt. $$
\noindent Again invoking Lemma \ref{lemfbound}, it then follows that $C$ can be chosen larger if necessary, still independent of $d_u, k>0$, so that $d_u \int_0^T \int_{- \infty}^0 (u^k_x)^2 \leq C/2$. If $d_v >0$, the estimate for $v_x^2$ can be proved likewise, using the equation for $v_k$. \qed

\medskip \bigskip

\noindent The following estimates for the differences of space and time translates of solutions  will yield sufficient compactness both  to obtain  the existence of solutions of $(P_1^k)$ when $d_u>0$ and $d_v=0$, and to study the strong-interaction limit ($k \to \infty$) of $(P_1^k)$. 
Here we want to allow $d_v=0$ and do not have $L^2(Q_T)$ bounds for $\vk_x$ in this case. Thus we cannot simply refer to \cite[Lemma 2.6]{CDHMN} to control the differences of space translates, as in the proof of Lemma \ref{lemspace}, but instead need an alternative method. 
Our proof centres on showing that  solutions $(\uk, \vk)$ of $(P_1^k)$ satisfy the $L^1$-contraction property (\ref{est1}). Note that the monotonicity properties of $F$ are used here, in establishing the sign condition (\ref{nondec}). See also \cite[Prop. 4]{hkr} and \cite[Prop. 3]{hmm} for some related arguments.

\bigskip
\noindent It is convenient to introduce a shorthand notation for space and time translates. Given a function $h$, let
\begin{equation}
\label{spacetimedef}
S_{\xi} h(x,t) : = h(x+\xi, t), \;\;\; T_{\tau} h(x,t) : = h(x, t+\tau),
\end{equation}

\noindent for all $(x,t)$ in a suitable space-time domain and appropriate $\xi$ and $\tau$.

\medskip

\begin{lemma}
\label{leml1contraction}
Suppose that $d_u > 0$ and $d_v \geq 0$, and let $(\uk, \vk)$  be a solution of $(P_1^k)$ satisfying (\ref{linfbound}). Then there exists a function $G \geq 0$, independent of $d_v \geq 0$ and $k >0$, such that $G(\xi) \to 0$ as $|\xi| \to 0$, and for all $t \in (0, T)$, 
\begin{equation}
\label{l1contraction}
\int_{\R} | \uk(x,t)  - S_{\xi} \uk (x,t)| + | \vk(x,t) - S_{\xi} \vk (x,t)| ~dx \, \leq \, G(\xi).
\end{equation}
\end{lemma}

\pf  Define 
\begin{equation}
\label{ushiftdef}
u := \uk - S_{\xi} \uk , \qquad v := \vk - S_{\xi} \vk , \qquad u_0: = \uk_0- S_{\xi} \uk_0 , \qquad v_0:= \vk_0 - S_{\xi} \vk_0,
\end{equation}
\noindent so that
$$
\begin{array}{rlrl}\displaystyle
u_t&= d_u u_{xx} -k\{ F(\uk, \vk) - F(S_{\xi}\uk, S_{\xi}\vk)\} \quad &\mbox{in}~\qt
\\[.1 in]
\displaystyle v_t&= d_v v_{xx} -k\{ F(\uk, \vk) - F(S_{\xi}\uk, S_{\xi}\vk)\} \quad &\mbox{in}~\qt\\[.1 in]
 u(x,0)&= \uk_0(x) - \uk_0(x + \xi) ,\quad v(x,0)=\vk_0(x) - \vk_0(x + \xi)   &\qquad\mbox{for} \quad x\in \R.
\end{array}$$

\medskip
\noindent Now given $L, \alpha>0$, let $\psi^L$ be the cut-off functions defined before Lemma \ref{lemkuvest}, and let $\ma$ be as defined in the proof of Lemma \ref{lemloneboundwholeline}. Then multiplying the equation for $u$ by $\psi^L \ma'(u)$ and integrating over $\R \times (0, t_0)$, $t_0 \in (0, T)$, gives
\begin{eqnarray}
\lefteqn{\int_{\R} \psi^L(x) \ma(u(x, t_0))~dx  =   \int_{\R} \psi^L(x) \ma(u_0(x))~dx - d_u \int_0^{t_0} \int_{\R} \psi^L \ma''(u) (u_x)^2~dxdt } \nonumber \\
& & \hspace*{2cm}+ d_u \int_0^{t_0} \int_{\R} \psi^L_{xx} \ma(u)~dxdt  - \int_0^{t_0} \int_{\R} \psi^L \ma'(u) k \{  F(\uk, \vk) - F(S_{\xi}\uk, S_{\xi}\vk)\}~dxdt \nonumber \\
& \leq & \int_{\R} \psi^L(x) \ma(u_0(x))~dx+ d_u \int_0^{t_0} \int_{\R} \psi^L_{xx} \ma(u)~dxdt \nonumber \\
& & \hspace*{5cm} - \int_0^{t_0} \int_{\R} \psi^L \ma'(u) k \{  F(\uk, \vk) - F(S_{\xi}\uk, S_{\xi}\vk) \}~dxdt,
\label{lastline}
\end{eqnarray}

\noindent \textcolor{black}{since $\ma''(u) \geq 0$ because $\ma$ is convex.} Then letting $\alpha \to 0$ in (\ref{lastline}) yields that for each $L>0$ and each $t_0 \in (0, T)$, 
\begin{eqnarray}
\lefteqn{\int_{\R} \psi^L(x) |u(x, t_0)|~dx \leq \int_{\R} \psi^L(x) |u_0(x)|~dx + d_u \int_0^{t_0}  \int_{\R} \psi^L_{xx}(x) |u(x,t)|~dxdt} \nonumber \\ & & \hspace*{4cm} - \int_0^{t_0} \int_{\R} \psi^L(x) \mbox{sgn}(u) k \{  F(\uk, \vk) - F(S_{\xi}\uk, S_{\xi}\vk)\}~dxdt, \label{uest}
\end{eqnarray}

\noindent and similarly,
 \begin{eqnarray}
\lefteqn{\int_{\R} \psi^L(x) |v(x, t_0)|~dx \leq \int_{\R} \psi^L(x) |v_0(x)|~dx + d_v \int_0^{t_0}  \int_{\R} \psi^L_{xx}(x) |v(x,t)|~dxdt} \nonumber \\ & & \hspace*{4cm} - \int_0^{t_0} \int_{\R} \psi^L(x) \mbox{sgn}(v) k \{  F(\uk, \vk) - F(S_{\xi}\uk, S_{\xi}\vk)\}~dxdt. \label{vest}
\end{eqnarray}

\noindent Adding (\ref{uest}) and (\ref{vest}) then gives 
 \begin{eqnarray}
\lefteqn{\int_{\R} \psi^L(x) \{|u(x, t_0)|+ |v(x, t_0)|\}~dx  \leq  \int_{\R} \psi^L(x) \{ |u_0(x)| + |v_0(x)|\} ~dx} \nonumber \\ &  & \hspace*{4cm} + \int_0^{t_0}  \int_{\R} \psi^L_{xx}(x)\{ d_u |u(x,t)| + d_v |v(x,t)| \} ~dxdt \label{totalest} \\ &  & \hspace*{4cm} - k \int_0^{t_0} \int_{\R} \psi^L(x) \{ \mbox{sgn}(u) + \mbox{sgn}(v)\}  \{  F(\uk, \vk) - F(S_{\xi}\uk, S_{\xi}\vk) \}~dxdt. \nonumber
\end{eqnarray}

\noindent Then since $F(\cdot, v)$ and $F(u, \cdot)$ are non-decreasing for all $u,v \in \R^+$, we have
\begin{equation}
\label{nondec}
(\mbox{sgn}(u) + \mbox{sgn}(v) ) \{   F(\uk, \vk) - F(S_{\xi}\uk, S_{\xi}\vk) \} \geq 0,
\end{equation}

\noindent because either $\mbox{sgn}(u) + \mbox{sgn}(v)=0$ or else $\mbox{sgn}(u) + \mbox{sgn}(v)$ and $F(\uk, \vk) - F(S_{\xi}\uk, S_{\xi}\vk) $ have the same sign, and hence
\begin{eqnarray}
\label{nextest}
\lefteqn{\int_{\R} \psi^L(x) \{ |u(x, t_0)| + |v(x, t_0)| \}~dx \leq \int_{\R} \psi^L(x) \{ |u_0(x)| + | v_0(x)|\}~dx} \\ & & \hspace*{5cm} + \int_0^{t_0} \int_{\R} \psi^L_{xx}(x) \{ d_u|u(x,t)| + d_v |v(x,t)| \} ~dxdt. \nonumber 
\end{eqnarray}

\noindent Now by Lemma \ref{lemloneboundwholeline},  $u(\cdot, t) - u^{\infty}_0, v(\cdot, t)- v^{\infty}_0 \in L^1(\R)$ for each
$t \in (0, T)$. Hence we can  let $L \to \infty$ in (\ref{nextest}) and thus obtain that for each $t_0 \in (0, T)$, 
\begin{eqnarray}
\label{est1}
\lefteqn{\int_{\R}  \{ |\uk(x, t_0) - \uk(x + \xi, t_0)| + |\vk(x, t_0) - \vk(x + \xi, t_0)| \}~dx} \\ & &  \hspace*{4cm} \leq \int_{\R}  \{ |\uk_0(x) - \uk_0(x + \xi)| + | \vk_0(x) - \vk_0(x + \xi)|\}~dx. \nonumber
\end{eqnarray}
\noindent The existence of $G$ is then immediate  from the assumption that 
$\| \uk_0 (\cdot + \xi) - \uk_0 (\cdot) \|_{L^1(\R)} + \| \vk_0 (\cdot + \xi) - \vk_0 (\cdot) \|_{L^1(\R)} \; \leq \; \omega(|\xi|)$ where $\omega(|\xi|) \to 0$ as $\xi \to 0$.

   \qed

\bigskip
\begin{lemma}
\label{lemtimetranslates}
Suppose that $d_u > 0$, $d_v \geq 0$ and let $(u^k, v^k)$ be a solution of $(P_1^k)$
satisfying (\ref{linfbound}). Then there exists $C>0$, independent of  $d_v$ and $k$, such
that for any $\tau \in (0, T)$, 
\begin{eqnarray*}
\int_0^{T- \tau} \int_{\R} | T_{\tau}\uk(x,t)- u^k(x,t)|^2 ~ dx dt & \leq & \tau C,\\
\int_0^{T- \tau} \int_{\R} | T_{\tau}\vk(x,t) - v^k(x,t)|^2 ~ dx dt & \leq & \tau C.
\end{eqnarray*}
\end{lemma}

\pf The proof is similar to that of \cite[Lemma 3]{hkr}; see also \cite[Lemma 2.6]{CDHMN}.
We sketch the key points here, focussing on the parts where our problem needs a slightly different argument.
Let $\psi^L$ be as defined before Lemma \ref{lemkuvest}. 
Then it follows using the equation for $u^k$ that
\begin{eqnarray}
\lefteqn{\int_0^{T- \tau} \int_{\R}  [T_{\tau}\uk(x,t) - u^k(x,t)]^2\psi^L~dxdt} \nonumber \\
& = & \int_0^{\tau} \int_0^{T- \tau} \int_{\R} 
(u^k(x, t+ \tau) - u^k(x, t)) [ d_u u^k_{xx}(x, t+s) - k F(u^k, v^k)] \psi^L ~ dxdtds \nonumber \\
& = &- \int_0^{\tau} \int_0^{T- \tau} \int_{\R} d_u [u^k(x, t+ \tau) - u^k(x,t)] u_x(x, t+s) \psi^L_x~dxdtds
\nonumber \\
&  &  - \int_0^{\tau} \int_0^{T- \tau} \int_{\R} d_u [u^k(x, t+ \tau) - u^k(x,t)]_x u_x(x, t+s) \psi^L~dxdtds
\nonumber \\
&  & - \int_0^{\tau} \int_0^{T- \tau} \int_{\R} [u^k(x, t+ \tau) - u^k(x,t)] k F(u^k, v^k)\, \psi^L ~ dxdtds 
\nonumber \\
& \leq & \left(\sup | \psi^L_x| \right)  \int_0^{\tau} \int_0^{T- \tau} \int_{L \leq |x| \leq L+1} d_u |u^k(x, t+ \tau) - u^k(x,t)| |u_x(x, t+s)|~dxdtds
\nonumber \\ & &
\hspace*{1cm}+\; 2 d_u \tau \int_0^T \int_{\R} (u^k_x)^2(x,t) ~ dx dt 
\; + \; 2M \tau \int_0^T \int_{\R} k F(u^k, v^k)~dxdt. \label{newpart}
\end{eqnarray}

\noindent Now the mapping $(x,t) \mapsto |u^k(x, t+ \tau) - u^k(x,t)| |u_x(x, t+s)|$ is integrable on $\R \times (0, T - \tau)$, by (\ref{linfbound}), Lemma \ref{lemloneboundwholeline}  and Lemma \ref{leml2derivbound}, and $\left(\sup | \psi^L_x| \right)$ is bounded independently of $L$. So the first term on the right-hand side of (\ref{newpart}) tends to $0$ as $L \to \infty$. Thus letting
$L \to \infty$ yields
$$\int_0^{T- \tau} \int_{\R}  [u^k(x, t+ \tau) - u^k(x,t)]^2~dxdt \leq 2 d_u \tau \int_0^T \int_{\R} (u^k_x)^2(x,t) ~ dx dt 
\; + \; 2M \tau \int_0^T \int_{\R} k F(u^k, v^k)~dxdt,$$
\noindent from which the estimate for $u^k$ follows using Lemma \ref{lemfbound} and Lemma \ref{leml2derivbound}.  When $d_v>0$, the estimate for $v^k$ follows likewise, using the equation for $v^k$. When $d_v =0$, a similar but simpler argument applies, omitting the terms deriving from   $v^k_{xx}$.\qed

\bigskip \bigskip
\noindent We can now prove a convergence result for solutions $(u^k, v^k)$ of $(P_1^k)$ as $d_v \to 0$. 

\begin{lemma}
\label{lemconvergencedzero} Let $k>0$ and $d_u>0$ be fixed and $(u^k_{d_v}, v^k_{d_v})$ be solutions of $(P_1^k)$ satisfying (\ref{linfbound}) with $d_v >0$. Then there exists $(u^k_*, v^k_*) \in (L^{\infty}(Q_T))^2$ such that up to a subsequence, for each $J>0$, 
$$ \begin{array}{cl}
u^k_{d_v} \to u^k_* & \mbox{in} \;\; L^2((-J, J) \times (0, T)),\\
v^k_{d_v} \to v^k_* & \mbox{in} \;\; L^2((-J, J) \times (0, T)),\\
u^k_{d_v} - \tilde{u} \rightharpoonup u^k_* - \tilde{u} & \mbox{in} \;\; L^2(0, T; H^1(\R)),
\end{array}
$$
\noindent as $d_v \to 0$, where $\tilde{u} \in C^{\infty}(\R)$ is a smooth function such that $\tilde{u}(x) = u_0^{\infty}(x)$ for all $|x| \geq 1$. 
\end{lemma}

\pf It follows from Lemma \ref{lemloneboundwholeline} and (\ref{linfbound}) that $\|u^k_{d_v} - u^{\infty}_0\|_{L^2(Q_T)}$ and $\|v^k_{d_v} - v^{\infty}_0\|_{L^2(Q_T)}$ are  bounded independently of $d_v$. So  Lemmas \ref{leml1contraction}, \ref{lemtimetranslates} and the Riesz-Fr\'echet-Kolmogorov Theorem \cite[Theorem 4.26]{brezisnew} yield that the sets $\{ u^k_{d_v} - u^{\infty}_0\}_{d_v>0}$ and
$\{ v^k_{d_v} - v^{\infty}_0\}_{d_v>0}$ are each relatively compact in $L^2((-J, J) \times (0, T))$ for each $J>0$. The weak convergence of $u^k_{d_v} - \tilde{u}$ in $L^2(0, T; H^1(\R))$ follows from  the fact that $\|u^k_{d_v} - u^{\infty}_0\|_{L^2(Q_T)}$ is bounded independently of $d_v$ together with Lemma \ref{leml2derivbound}. \qed

\bigskip

\begin{theorem}
\label{thmexistzero} Let $d_v=0$ and $k>0$. Then problem $(P_1^k)$ has a unique weak solution
$$(\uk, \vk) \in W^{2,1}_p( (-J, J) \times (0, T)) \times W^{1,\infty}(0, T; L^{\infty}((-J, J)))\;\; \mbox{for each} \; J>0 , \; \; p \geq 1,$$
where $(\uk, \vk)$ is a weak solution in the sense that 
\begin{eqnarray}
\iint_{Q_T} \uk \psi_t ~dxdt  + \iint_{Q_T} \{ d_u \uk \psi_{xx} - k F(\uk, \vk) \psi \}~dxdt & = & - \int_{\R} u_0^k
\psi( \cdot, 0) ~ dx, \label{udz}\\
\iint_{Q_T} \vk \psi_t ~dxdt  - \iint_{Q_T} k F(\uk, \vk) \psi ~dxdt & = & - \int_{\R} v_0^k
\psi( \cdot, 0) ~ dx, \label{vdz}
\end{eqnarray}
for all $\psi \in \mathcal{F}_T = \{\psi \in C^{2,1}(Q_T): \psi (\cdot, T) =0\; \mbox{and} \;  \mbox{supp}\, \psi \subset [-J, J] \times [0, T] \;\mbox{for some} \; J>0 \}$, and also satisfies  $0 \leq \uk, \vk \leq M$.
\end{theorem}

\medskip
\pf Multiplying $(P_1^k)$ by $\psi \in \mathcal{F}_T$ and integrating over $Q_T$ yields that for each $d_v >0$, solutions $(u^k_{d_v}, v^k_{d_v})$ of  $(P_1^k)$ satisfy 
\begin{eqnarray}
\iint_{Q_T} u^k_{d_v} \psi_t ~dxdt  + \iint_{Q_T} \{ d_u u^k_{d_v} \psi_{xx} - k F(u^k_{d_v}, v^k_{d_v}) \psi \}~dxdt & = & - \int_{\R} u_0^k
\psi( \cdot, 0) ~ dx, \label{udz1}\\
\iint_{Q_T} v^k_{d_v} \psi_t ~dxdt  + \iint_{Q_T} \{  d_v v^k_{d_v} \psi_{xx}  - k F(u^k_{d_v}, v^k_{d_v}) \psi \} ~dxdt & = & - \int_{\R} v_0^k
\psi( \cdot, 0) ~ dx. \label{vdz1}
\end{eqnarray}
Then the existence of a solution $(\uk, \vk)$ to (\ref{udz})-(\ref{vdz}) follows by using Lemma \ref{lemconvergencedzero}
to pass to the limit along a subsequence as $d_v \to 0$ in (\ref{udz1})-(\ref{vdz1}). The regularity  of $\uk$ follows from the fact that solutions $(u^k_{d_v}, v^k_{d_v})$ of $(P_1^k)$ satisfying  (\ref{linfbound}) are such that $(u^k_{d_v})_t - d_u (u^k_{d_v})_{xx} = - k F(u^k_{d_v}, v^k_{d_v})$ is bounded in $L^{\infty}(Q_T)$ independently of $d_v$, which, since $u_0^k \in C^2(\R)$, implies that 
$u^k_{d_v}$ is bounded independently of $d_v>0$ in $W^{2,1}_p((-J, J) \times (0, T))$ for each $J>0$ and $p \geq 1$. The regularity of $\vk$ is immediate from the fact that (\ref{vdz}) implies that $\vk \in W^{1, \infty}(0, T; L^{\infty}((-J, J)))$   for each $J>0$, and the uniqueness of $(\uk, \vk)$ follows from the comparison principle in  Lemma \ref{lemcomparisonmain}. \qed

\label{subsectiondvzerowhole}

\bigskip \medskip

\bigskip 

\subsection{The limit problem for $(P_1^k)$ as $k \to \infty$}
\label{subsectionkinftywhole}

The a priori estimates of the previous section yield sufficient compactness to establish the existence of limits of solutions of $(P_1^k)$ as $ k \to \infty$, both when $d_v>0$ and when $d_v=0$. \textcolor{black}{The proof of the following result is directly analogous to that of Lemma \ref{lemconvergencedzero}, using bounds independent of $k$ in place of bounds independent of $d_v$, and is left to the reader.}
\smallskip

\begin{lemma}
\label{lemkinftyconvergence}
Let $d_u >0$ and $d_v \geq 0$ be fixed and $(u^k, v^k)$ be solutions of $(P_1^k)$ satisfying (\ref{linfbound}) with $k>0$. Then there exists $(u, v) \in (L^{\infty}(Q_T))^2$ such that up to a subsequence, for each $J>0$, 
$$ \begin{array}{cl}
u^k \to u & \mbox{in} \;\; L^2((-J, J) \times (0, T)),\\
v^k \to v & \mbox{in} \;\; L^2((-J, J) \times (0, T)),\\
u^k - \tilde{u} \rightharpoonup u - \tilde{u} & \mbox{in} \;\; L^2(0, T; H^1(\R)),
\end{array}
$$
\noindent as $k \to \infty$, where $\tilde{u} \in C^{\infty}(\R)$ is a smooth function such that $\tilde{u}(x) = u_0^{\infty}(x)$ for all $|x| \geq 1$. 
\end{lemma}

\medskip
\noindent The following segregation result is a key to the characterisation of the limits $u, v$ in Lemma \ref{lemkinftyconvergence}.

\begin{lemma}
\label{lemsegregation}
Let $d_u >0$, $d_v \geq 0$ and $(u,v)$ be as in Lemma \ref{lemkinftyconvergence}. Then
\begin{equation}
\label{seg}
uv =0 \;\;\; a.e. \; \mbox{in} \;\; Q_T.
\end{equation}
\end{lemma}

\pf It follows from Lemmas \ref{lemfbound} and \ref{lemkinftyconvergence} that $F(u,v)=0$ almost everywhere in $Q_T$, from which (\ref{seg}) follows since $F(u, v)=0$ if and only if $u=0$ or $v=0$. \qed

\bigskip
\noindent To derive the limit problem, set 
\begin{equation}
\label{wkwdef}
\wk : = \uk - \vk, \;\;\;\; w : = u-v.
\end{equation}

\noindent Then it follows from Lemmas \ref{lemkinftyconvergence} and \ref{lemsegregation} that as a sequence $k_n \to \infty$, 
$$w^{k_n} \to w \;\; \mbox{in} \;\; L^2((-J, J) \times (0, T)) \;\; \mbox{for all} : J>0\;\; \mbox{and} \; a.e. \; \mbox{in} \; Q_T, $$
and that
$$ u=w^+ \;\; \mbox{and} \;\; v = -w^-,$$
where $s^+ = \max \{0, s\}$ and $s^- = \min \{0, s\}.$ Next note the following equality. 

\bigskip

\begin{lemma}
\label{lemlimiteq}
Let $d_u >0$, $d_v \geq 0$ and $(u,v)$ be as in Lemma \ref{lemkinftyconvergence}. Then
\begin{equation}
\label{limiteq}
- \iint_{Q_T} (u-v) \psi_t~dxdt  - \int_{\R} (u^{\infty}_0 - v^{\infty}_0) \, \psi(x, 0)~dx = \iint_{Q_T}(d_u u - d_v v) \
\psi_{xx} ~dxdt,
\end{equation}
for all $\psi \in \mathcal{F}_T = \{\psi \in C^{2,1}(Q_T): \psi (\cdot, T) =0\; \mbox{and} \;  \mbox{supp}\, \psi \subset [-J, J] \times [0, T] \;\mbox{for some} \; J>0 \}$. 

\end{lemma}

\pf Multiplying the difference between the equations for $\uk$ and $\vk$ by $\psi \in \mathcal{F}_T$ and integrating over $Q_T$ yields
$$- \iint_{Q_T} (\uk-\vk) \psi_t~dxdt  - \int_{\R} (u^k_0 - v^k_0) \, \psi(x, 0)~dx = \iint_{Q_T}(d_u \uk - d_v \vk) \
\psi_{xx} ~dxdt,$$
from which (\ref{limiteq}) follows using Lemma \ref{lemkinftyconvergence} and the fact that $\uk_0 \to u_0^{\infty}$
and $\vk_0 \to v_0^{\infty}$ in $L^1(\R)$ as $k \to \infty$. \qed

\bigskip
\noindent Now define
\begin{equation}
\label{defcurld}
\mathcal{D}(s) : = \left\{ \begin{array}{ll} d_u s & \mbox{if} \; s \geq 0,\\ d_v s & \mbox{if} \; s <0,   \end{array}\right.
\end{equation}
and the limit problem
$$(P_1^{limit}) \left\{\begin{array}{rlrl}\displaystyle
w_t&= \mathcal{D}(w)_{xx},\quad &\mbox{in}~\R \times [0, \infty),
 \\[.1 in]
 w(x,0)&= w_0(x) : = \left\{ \begin{array}{c} U_0 \\ - V_0,  \end{array} \right. & \begin{array}{l}
 \mbox{if}\; x<0, \\ \mbox{if} \; x>0.\end{array} 
 \end{array} \right.$$
 
 \bigskip
 
 \begin{definition}
 \label{defweaksolution}
 A function $w$ is a weak solution of Problem $(P_1^{limit})$ if
 \begin{itemize}
 \item[(i)] $w \in L^{\infty} (\R \times \R^+)$,
 \item[(ii)] for all $T>0$, 
 $$\iint_{Q_T} (w \psi_t + \mathcal{D}(w) \psi_{xx})~dxdt = - \int_{\R} w_0(x) \psi(x,0)~dx,$$
 for all $\psi \in \mathcal{F}_T $ $$ = \{\psi \in C^{2,1}(Q_T): \psi (\cdot, T) =0\; \mbox{and} \;  \mbox{supp}\, \psi \subset [-J, J] \times [0, T] \;\mbox{for some} \; J>0 \}.$$ 
  \end{itemize}
 \end{definition}
 
 \bigskip

 \begin{lemma}
 \label{lemuniqueweaklimit}
 The function $w$ defined in equation (\ref{wkwdef}) is the unique weak solution of Problem $(P_1^{limit})$ and the whole sequence $(\uk, \vk)$ in Lemma
 \ref{lemkinftyconvergence} converges to $(w^+, -w^-)$.
 \end{lemma}
 
 \pf That $w$ is a weak solution of $(P_1^{limit})$ follows immediately from Lemma \ref{lemlimiteq} and the definition of $\mathcal{D}$. The uniqueness is a consequence of
 \cite[Appendix, Proposition A]{bkp}, which extends the method of \cite[Proposition 9]{acp} to unbounded domains, via exactly the argument used to establish uniqueness for a  similar problem in \cite[Appendix, Proof of Theorem C]{bkp}.  Note that although it is assumed throughout \cite[Appendix]{bkp} that the initial data of the problems considered is continuous, it is straightforward to verify that this is not in fact necessary for the proofs. \qed

\bigskip 
\noindent We next identify the limit $w$ as a certain self-similar solution of Problem $(P_1^{limit})$, the precise form of which depends on whether $d_v >0$ or
$d_v=0$. To this end, we first state a free-boundary problem, including interface conditions, that is satisfied by the solution $w$ of $(P_1^{limit})$ under 
some regularity assumptions and conditions on the form of the free boundary. The proof follows immediately from that of   \cite[Theorem 5]{hmm}
and we omit it.

\begin{theorem}
\label{thmfreeboundary}
Let $w$ be the unique weak solution of Problem $(P_1^{limit})$. Suppose that there exists a function $\xi: [0, T] \to \R$ such that  for each $t \in [0, T]$, 
$$ w(x, t) > 0 \;\;\; \mbox{if} \; \;\; x < \xi(t) \;\;\; \mbox{and} \;\;\; w(x, t) < 0 \;\;\; \mbox{if} \;\;\; x> \xi(t). $$
Then if $t \mapsto \xi(t)$ is sufficiently smooth and the functions $u := w^+$ and $v: = - w^-$ are smooth up to $\xi(t)$, the functions
$u$ and $v$ satisfy 
$$ 
(P_1^{limit}) \left\{
\begin{array}{ll}
u_t = d_u u_{xx},  & \mbox{in} \;\; \{(x,t) \in Q_T: x <  \xi(t)\} ,\\
v_t = d_v v_{xx},   & \mbox{in} \;\; \{(x,t) \in Q_T: x >  \xi(t)\} ,\\
\left[u \right] = d_v \left[v \right] =0,     & \mbox{on} \;\; \Gamma_T  := \{(x,t) \in Q_T: x =  \xi(t)\},\\
\left [v \right] \xi'(t) = \left[ d_u u_x - d_v v_x \right], & \mbox{on} \;\; \Gamma_T := \{(x,t) \in Q_T: x =  \xi(t)\},\\
u(\cdot, 0) = u_0^{\infty} (\cdot), & \mbox{in} \;\; \R,\\
v(\cdot, 0) = v_0^{\infty}(\cdot),  & \mbox{in} \;\; \R,
\end{array}
\right.
$$

\noindent where $\left[ \cdot \right]$ denotes the jump across $\xi(t)$ from $\{x < \xi(t) \}$ to $\{x > \xi(t)\}$, that is, $\left[a\right] : = \lim_{x \downarrow \xi(t)} a(x, t) - \lim_{x \uparrow \xi(t)} a(x,t)$, $\xi'(t)$ denotes the speed of propagation of the free boundary $\xi (t)$.
\end{theorem}

\bigskip

\noindent Interpreting the interface conditions on $\Gamma_T$ then yields the following two limit problems. 

\begin{corollary}
\label{corfreeboundary}
Let $w$ and $\xi :[0,T] \to \R$ satisfy the hypotheses of Theorem \ref{thmfreeboundary}. 
Then  the functions
$u:=w^+$ and $v=-w^-$ satisfy one of two limit problems, depending on whether $d_v>0$ or $d_v=0$. If $d_v>0$, then 
$$ 
(P^{limit}_{1, d_v>0}) \left\{
\begin{array}{ll}
u_t = d_u u_{xx},  & \mbox{in} \;\; \{(x,t) \in Q_T: x <  \xi(t)\} ,\\
v = 0,  & \mbox{in} \;\; \{(x,t) \in Q_T: x <  \xi(t)\} ,\\
v_t = d_v v_{xx},   & \mbox{in} \;\; \{(x,t) \in Q_T: x >  \xi(t)\} ,\\
u = 0,   & \mbox{in} \;\; \{(x,t) \in Q_T: x >  \xi(t)\} , \\
\lim_{x \uparrow \xi(t)} u(x,t) =0 = \lim_{x \downarrow \xi(t)} v(x,t),    & \mbox{for each} \;\; t \in [0, T], \\
d_u \lim_{x \uparrow \xi(t)} u_x(x,t) = - d_v \lim_{x \downarrow \xi(t)} v_x(x,t), & \mbox{for each} \;\; t \in [0, T],\\
u(\cdot, 0) = u_0^{\infty} (\cdot), & \mbox{in} \;\; \R,\\
v(\cdot, 0) = v_0^{\infty}(\cdot),  & \mbox{in} \;\; \R,
\end{array}
\right.
$$
\noindent whereas if $d_v=0$ and  we suppose additionally that 
$\xi(0)=0$ and $t \mapsto \xi(t)$ is a non-decreasing function, then 
$$ 
(P^{limit}_{1, d_v=0}) \left\{
\begin{array}{ll}
u_t = d_u u_{xx},  & \mbox{in} \;\; \{(x,t) \in Q_T: x <  \xi(t)\} ,\\
v = 0,  & \mbox{in} \;\; \{(x,t) \in Q_T: x <  \xi(t)\} ,\\
v = V_0,   & \mbox{in} \;\; \{(x,t) \in Q_T: x >  \xi(t)\} ,\\
u = 0,   & \mbox{in} \;\; \{(x,t) \in Q_T: x >  \xi(t)\} ,\\
\lim_{x \uparrow \xi(t)} u(x,t) =0,  & \mbox{for each} \;\; t \in [0, T],   \\
V_0 \; \xi'(t) = - {d_u}\lim_{x \uparrow \xi(t)} u_x (x,t), & \mbox{for each} \;\; t \in [0, T],\\
u(\cdot, 0) = u_0^{\infty} (\cdot), & \mbox{in} \;\; \R,\\
v(\cdot, 0) = v_0^{\infty}(\cdot),  & \mbox{in} \;\; \R,
\end{array}
\right.
$$
where 
$\xi'(t)$ denotes the speed of propagation of the free boundary $\xi (t)$. 
\end{corollary}

\bigskip

\pf We interpret the meaning of the interface conditions in Theorem \ref{thmfreeboundary}, depending on whether $d_v>0$ or $d_v=0$. 
The  condition $$\left[u \right] = d_v \left[v \right] =0\;\;\; \mbox{on}\;\;\; \Gamma_T  := \{(x,t) \in Q_T: x =  \xi(t)\},$$
implies that
$u(\cdot, t)$ is continuous across $\xi(t)$, so that $$\lim_{x \uparrow \xi(t)} u(x,t) =  \lim_{x \downarrow \xi(t)} u(x,t)= 0 .$$
Moreover,  if $d_v>0$, then $v(\cdot, t)$ is also continuous across $\xi(t)$, and so $$\lim_{x \downarrow \xi(t)} v(x,t) = \lim_{x \uparrow \xi(t)} v(x,t) = 0 ,$$
whereas if $d_v=0$, $v(\cdot, t)$ may jump across $\xi(t)$. Indeed, since $\xi(0)=0$ and $t \mapsto \xi(t)$ is a non-decreasing function, 
it follows from the fact that $v_t=0$ in $\{(x,t) \in Q_T: x >  \xi(t)\} $  if $d_v=0$, together with the initial condition that $v_0(x) = V_0$ if $x >0$, that
$v(x,t) \equiv V_0$ for all $x \geq \xi(t)$, and thus
$$\left[ v \right] = V_0 - 0 = V_0\;\;\; \mbox{for all} \;\;t \in [0, T].$$

\noindent The normal derivative condition $$\left [v \right] \xi'(t) = \left[ d_u u_x - d_v v_x \right] \;\;\; \mbox{on} \;\; \Gamma_T := \{(x,t) \in Q_T: x =  \xi(t)\},$$
implies that  if $d_v>0$, then $0  = \left[ d_u u_x - d_v v_x \right],$
which says that
$$d_u \lim_{x \uparrow \xi(t)} u_x(x,t) = - d_v \lim_{x \downarrow \xi(t)} v_x(x,t), $$
or equivalently, $$d_u \lim_{x \uparrow \xi(t)} w^+_x(x,t) = 
d_v \lim_{x \downarrow \xi(t)} w^-_x(x,t).$$
\bigskip
\noindent 
On the other hand, if $d_v=0$, then 
$$\lim_{x \downarrow \xi(t)} v(x,t)\;\;  \xi'(t) = - d_u \lim_{x \uparrow \xi(t)} u_x (x,t),$$

\noindent which in the case that  $\xi(0)=0$ and $t \mapsto \xi(t)$ is a non-decreasing function gives
$$V_0 \; \xi'(t) = - {d_u}\lim_{x \uparrow \xi(t)} u_x (x,t).$$
\qed

\bigskip
\noindent It is then easy to show that the limit problems in Corollary \ref{corfreeboundary} admit self-similar solutions.

\begin{theorem}
\label{thmselfsimilar}
The unique weak solution $w$ of Problem $(P_1^{limit})$ has a self-similar form. There exists a function $f: \R \to \R$ and a constant $a \in \R$ such that
\begin{equation}
\label{wform}
w(x,t) = f\left( \frac{x}{\sqrt{t}}\right) , \;\;\; (x,t) \in Q_T, \;\;\; \mbox{and} \;\;\;\xi(t) = a \sqrt{t}, \;\;\; t \in [0, T].
\end{equation}
\noindent If $d_v>0$, then $a \in \R$ is the unique root of the equation
$$  d_u U_0 \int_a^{\infty} e^{\frac{a^2 - s^2}{4d_v}} ds = d_v V_0 \int_{- \infty}^{a} e^{\frac{a^2 - s^2}{4d_u}} ds,$$
\noindent and
\begin{equation}
\label{fone}
f(\eta) = \left\{ \begin{array}{ll} U_0 \left(  1 - \frac{\int_{- \infty}^{\eta} e^{- \frac{s^2}{4d_u}}\;ds}{\int_{- \infty}^{a} e^{- \frac{s^2}{4d_u}} \;ds}        \right), & \mbox{if} \;\; \eta \leq a,\\
-V_0 \left(  1 - \frac{\int_{\eta}^{\infty} e^{- \frac{s^2}{4d_v}}\;ds}{\int_{a}^{\infty} e^{- \frac{s^2}{4d_v}} \;ds}        \right), & \mbox{if} \;\; \eta>a.
 \end{array} \right.
 \end{equation}
\noindent On the other hand, if $d_v=0$, then $a>0$ is the unique root of the equation
$$U_0 = \frac{V_0 a}{2d_u} \int_{- \infty}^a e^{\frac{a^2 - s^2}{4 d_u}}\;ds,$$
\noindent and
\begin{equation}
\label{ftwo}
f(\eta) = \left\{ \begin{array}{ll} U_0 \left(  1 - \frac{\int_{- \infty}^{\eta} e^{- \frac{s^2}{4d_u}}\;ds}{\int_{- \infty}^{a} e^{- \frac{s^2}{4d_u}} \;ds}        \right), & \mbox{if} \;\; \eta \leq a,\\
-V_0 , & \mbox{if} \;\; \eta>a.
 \end{array} \right.
 \end{equation}
\end{theorem}

\bigskip

\pf Straightforward verification shows that the functions $w$ defined in (\ref{wform}) and (\ref{fone}) or (\ref{ftwo}) satisfy $(P^{limit}_{1, d_v>0})$ or $(P^{limit}_{1, d_v=0})$
when $d_v>0$ or $d_v=0$ respectively, and hence give a solution of the problem $(P^{limit})$, which must therefore be the unique solution. \qed

\bigskip
\noindent As already mentioned in the Introduction,  the constant $a$ is not necessarily positive in the case when $d_v>0$. Some sufficient conditions ensuring the sign of $a$ are as follows.

\medskip
\begin{proposition}
\label{propsuffpos}
Suppose that $d_u, d_v, U_0, V_0 \in \R$ are all strictly positive, and let $a \in \R$ be the unique root of the equation
\begin{equation}
\label{aidentity} d_u U_0 \int_a^{\infty} e^{\frac{a^2 - s^2}{4d_v}} ds = d_v V_0 \int_{- \infty}^{a} e^{\frac{a^2 - s^2}{4d_u}} ds.
\end{equation}
Then
\begin{itemize}
\item[(i)] if $d_u=d_v$ and $U_0=V_0$, then $a=0$;
\item[(ii)] if $d_u \leq d_v$ and $\sqrt{d_u}U_0 \leq \sqrt{d_v} V_0$, then $a<0$;
\item[(iii)] if $d_u \geq d_v$ and $\sqrt{d_u}U_0 \geq \sqrt{d_v} V_0$, then $a>0$.
\end{itemize}
\end{proposition}

\bigskip
\pf  A straightforward rearrangement of (\ref{aidentity}) gives
$$\frac{\sqrt{d_u}\, U_0}{\sqrt{d_v}\, V_0} \; e^{\frac{a^2}{2} \left( \frac{1}{d_v} - \frac{1}{d_u} \right)}
\int_{\frac{a}{2 \sqrt{d_v}}}^{\infty} e^{- t^2}~dt \; = \; \int_{- \infty}^{\frac{a}{2 \sqrt{d_u}}} e^{-t^2}~dt,$$
from which the result is  clear. \qed

\label{sec-klimitwhole}

\bigskip
\section{The half-line case: problem $(P_2^k)$}

\subsection{Existence and uniqueness of solutions of $(P_2^k)$ when $d_u>0$ and $d_v>0$}

\noindent Suppose that $d_u >0$ and $d_v >0$. Similarly to the whole-line case, we can use an approximate problem to establish existence of solutions of
$(P_2^k)$. Choose  $M \geq \max\{ U_0, V_0 \}$ and for each $R>1$, let $(P_2^{R, \mu})$ denote the problem 

$$
(P_2^{R, \mu}) \left\{\begin{array}{rlrl}\displaystyle
u_t&= d_u u_{xx} -kF_{\mu}(u,v)\quad &\mbox{in}~(0, R) \times (0, T),
\\
\displaystyle v_t&= d_v v_{xx} -kF_{\mu} (u,v)\quad &\mbox{in}~(0, R) \times (0,T), \\
u(0,t) & = U_0 & \mbox{for} \quad t \in (0, T),\\
u_x(R,t) & =0 & \mbox{for} \quad t \in (0, T),
\\
v_x(0,t) & = 0 & \mbox{for} \quad t \in (0, T),
\\
v_x(R,t) & =0 & \mbox{for} \quad t \in (0, T),
\\
 u(x,0)&=u^k_{0,R}(x),\quad v(x,0)=v^k_{0, R}(x)   &\qquad\mbox{for} \quad x\in (0, R),
\end{array}\right.$$

\noindent where $u^k_{0,R}, v^k_{0, R} \in C^2(\R^+)$ are such that $0 \leq u^k_{0,R} \leq M$, 
$0 \leq v^k_{0, R} \leq M$  and
\begin{eqnarray}
& & u^k_{0,R} (x) =0  \; \mbox{for}   \; x > \left(1 - {\frac{1}{R}} \right) R, \;\;\;\;\;\;\; v^k_{0,R} (x) =V_0  \; \mbox{for}   \; x > \left(1 - {\frac{1}{R}} \right) R,
   \label{uvrini}
   \end{eqnarray}
   
   \noindent which defines the functions $u^k_{0, R}$, $v^k_{0, R}$ on the half-line $(0, \infty)$. The regularisation $F_{\mu}$
   is as defined in section \ref{existuniquewhole}. 
   
   \medskip
   \noindent Arguments  analogous to those used in Section \ref{existuniquewhole} yield existence of solutions to problem $(P_2^k)$  by passing to the  limits $R \to \infty$ and $\mu \to 0$ in  problem $(P_2^{R, \mu})$. We omit repetition of the  details the proof and simply state the result. 
   
\begin{theorem}
\label{lemhalflineexistence}
Suppose that $d_u>0$ and $d_v>0$. Then given $k>0$, there exists a  classical solution $(u^k, v^k)$ of $(P_2^k)$, such that 
for each $\delta >0$, $J>0$ and $p \geq 1$, 
\begin{equation}
\label{reg-half}
\uk,\vk \in C^{2,1}(\R^+ \times [\delta, T] ) \cap C^0(\R^+ \times [0, T]) \cap W^{2,1}_p((0, J) \times (0, T)),
\end{equation}
and
\begin{equation}
\label{linfboundhalf}
0 \leq u^k \leq M, \;\;\; 0 \leq v^k \leq M \;\;\; \mbox{on} \;\;\; \R^+ \times (0, T).
\end{equation}
\end{theorem}

\medskip
\noindent Uniqueness is again a consequence of  a comparison theorem.

\begin{lemma}
\label{lemcomparisonhalf}
Suppose that $d_u > 0$, $d_v \geq 0$, and let $(\ou, \ov)$, $(\uu, \uv)$ be such that for each $J>0$ and $p \geq 1$, 
$\ou, \uu \in  L^{\infty}(S_T) \cap W^{2,1}_p((0, J) \times (0, T))$, $\ov, \uv \in L^{\infty}(S_T) \cap W^{2,1}_p((0, J) \times (0, T))$ if $d_v >0$, $\ov, \uv \in L^{\infty}(S_T) \cap W^{1, \infty}(0, T; L^{\infty}((-J, J)))$ if $d_v=0$, and $(\ou, \ov)$, $(\uu, \uv)$ satisfy
$$\begin{array}{l} \ou_t \geq d_u \ou_{xx} - kF(\ou, \ov),\\ \ov_t \leq d_v \ov_{xx} - kF(\ou, \ov),\end{array} \qquad 
\begin{array}{l} \uu_t \leq d_u \uu_{xx} - kF(\uu, \uv),\\ \uv_t \geq d_v \uv_{xx} - kF(\uu, \uv),\end{array} \quad \mbox{in} \;\; S_T, 
$$

$$ \ou (0, \cdot ) \geq \uu(0, \cdot), \qquad d_v \ov_x(0, \cdot) \geq d_v \uv_x(0,\cdot) \;\;\mbox{on}\;\; (0, T],   $$
and
$$\ou(\cdot, 0) \geq \uu(\cdot,0), \qquad \ov(\cdot,0) \leq \uv(\cdot, 0) \;\; \mbox{on} \;\;  \R^+.$$

\noindent Then 
$$
\ou \geq \uu \quad \mbox{and} \quad \ov \leq \uv \quad \mbox{in} \;\; S_T.
$$
\end{lemma}
\medskip

\pf This follows from the same form of argument used to show Lemma \ref{lemcomparisonmain}, replacing $Q_T$ with $S_T$, integrals over $\R$ with integrals
over $\R^+$, and the cut-off function $\psi^L$ by $\psi^L_+: = \psi^L|_{\R^+}$. We omit most of the details and only note two key calculations involving
the boundary $\{0\} \times (0, T)$. Taking $u:= \uu - \ou$ and $v: = \uv - \ov$, we have
$$u(0, \cdot) \leq 0 \;\;\; \mbox{and} \;\;\; d_v v_x(0, \cdot) \leq 0 \;\; \mbox{on} \;\; (0,T].$$
\noindent Thus $(m_{\alpha}^+)'(u(0,\cdot)) =0$, so that integrating over $\R^+ \times (0, t_0)$, $t_0 \in (0, T]$, gives
\begin{eqnarray*}
\lefteqn{\int_0^{t_0} \int_{\R^+} u_{xx} (m_{\alpha}^+)'(u) \psi^L_+ \; dxdt =} \\ & & -  \int_0^{t_0}\int_{\R^+} u_x \left[ (\ma^+)''(u) u_x \psi^L_+ + (\ma^+)'(u) \left(\psi^L_+\right)_x \right]~dxdt \, \leq \, \int_0^{t_0}\int_{\R^+} m_{\alpha}^+(u) \left(\psi^L_+\right)_{xx}~dxdt,
\end{eqnarray*}
\noindent whereas if $d_v >0$, then
\begin{eqnarray*}
\lefteqn{\int_0^{t_0}\int_{\R^+} v_{xx} (m_{\alpha}^+)'(v) \psi^L_+ \, dxdt} \\  & = & \int_0^{t_0} v_x (0, t) (m_{\alpha}^+)'(v(0,t))\;dt  -  \int_0^{t_0} \int_{\R^+} v_x \left[ (\ma^+)''(v) v_x \psi^L_+ + (\ma^+)'(v) \left(\psi^L_+\right)_x \right]~dxdt\\
& = & \int_0^{t_0} v_x (0, t) (m_{\alpha}^+)'(v(0,t)) \; dt - \int_0^{t_0}\int_{\R^+} (\ma^+)''(v) (v_x )^2 \psi^L_+ ~dxdt + \int_0^{t_0}\int_{\R^+} m_{\alpha}^+(v) \left( \psi^L_+ \right)_{xx}~dxdt \\
& \leq & \int_0^{t_0}\int_{\R^+} m_{\alpha}^+(v) \left( \psi^L_+ \right)_{xx}~dxdt.
\end{eqnarray*}
\noindent We refer the reader to  Lemma \ref{lemcomparisonmain} for the remainder of the proof. \qed

\bigskip
\begin{corollary}
\label{cor-half}
Suppose $d_u > 0$ and $d_v > 0$. Then given $k>0$, there is at most one solution $(u^k, v^k)$ of $(P_2^k)$ in $L^{\infty}(S_T) \cap W^{2,1}_p((0, J) \times (0, T))$ for each $J>0$, $p \geq 1$. 
\end{corollary}

\bigskip
\subsection{Existence and uniqueness of solutions for $(P_2^k)$ when $d_u>0$ and $d_v=0$}
\noindent Again we begin with some preliminary estimates, counterparts of results in section \ref{existuniquenodiffusion}. 
Here some different arguments are needed because of the boundary condition at $x=0$. 

%\begin{lemma}
%\label{lemuliboundhalfline}
%There exists a constant $C>0$, dependent on $d_u>0$ but independent of $d_v \geq 0$ and $k>0$, such that for any  solution $(\uk, \vk)$ of $(P_2^k)$ satisfying %(\ref{linfboundhalf}), we have
%\begin{equation}
%\label{heatbound}
%\int_{0}^{\infty} \uk (x,t)~dx \leq C\;\;\; \mbox{for all} \;\; t \in [0, T].
%\end{equation}
%\end{lemma}

%\pf Let $\overline{u}^k$ denote the solution of the heat equation
%$$
 %\left\{\begin{array}{rlrl}\displaystyle
%u_t&= d_u u_{xx} \quad &\mbox{in}~\st,
%\\[.1 in]
%\displaystyle u(0, t) &= U_0 \quad &\mbox{for all}~t \in [0, T], \\[.1 in]
 %u(x,0)&=u^k_0(x),   &\qquad\mbox{for} \quad x\in \R^+.
%\end{array}\right.$$

%\noindent Since $kF(\uk, \vk) \geq 0$,  the usual scalar comparison theorem  implies that $\uk \leq \overline{u}^k$
%on $S_T$. Then (\ref{heatbound}) follows from properties of solutions of the heat equation and the fact that $\|u_0^k\|_{L^1(\R^+)}$ is bounded independently of $k$.
%\qed

\medskip

\noindent The following key bound is the half-line counterpart of Lemma \ref{lemfbound}.

\begin{lemma}
\label{lemfboundhalfline}
There exists a constant $C>0$,  independent of  $d_v \geq 0$ and $k>0$,  such that for  any solution $(u^k, v^k)$ of $(P_2^k)$
satisfying (\ref{linfboundhalf}), we have
$$
\int_0^T \int_{0}^{\infty} k F(u^k, v^k) ~ dx\,dt \leq C.
$$
\end{lemma}

\pf  \textcolor{black}{Define a cut-off function $\beta \in C^{\infty}(\R^+)$ such that $0 \leq \beta(x) \leq 1$ for all $x \in \R^+$, $\beta^L(0) = \beta^L_x(0)=0$, $\beta(x) = 1$ for all $x \in [1,2]$, and $\beta(x) =0$ for $x \geq 3$.
Then given $L \geq 2$, define the family of cut-off functions $\beta^L \in C^{\infty}(\R)$ by $\beta^{L}(x) = \beta(x)$ when $x \in [0,1]$, $\beta^L(x) = 1$ when $x \in [1, L]$, and $\beta^L(x) = \beta(x-L+2)$ when $x \geq L$. Note  that $0 \leq \beta^L \leq 1$ for all $L$,  and 
$\beta^L_x$, $\beta^L_{xx}$ are bounded in both $L^{\infty}(\R^+)$ and $L^1(\R^+)$ independently of $L$. Also let $\hat{\beta} \in C^{\infty}(\R^+)$ be such that $0 \leq \hat{\beta} \leq 1$, $\hat{\beta}(x)=1$ for all $x \in [0, 1]$ and $\hat{\beta}(x)=0$ for all $x \geq 2$. }

\medskip
\noindent Then multiplying the equation for $\uk$ by $\beta^L$ and integrating over $\R^+ \times (0, t_0)$, $t_0 \in (0, T]$, gives that
\begin{eqnarray} \nonumber \lefteqn{\int_{\R^+} \beta^L(x)  u^k(x,t_0)\,dx +  \int_0^{t_0} \int_{\R^+} \beta^L(x) k F(u^k, v^k) \, dxdt  =} \\ & &  \hspace*{5cm}d_u \int_0^{t_0} 
\int_{\R^+} \beta^L_{xx}(x) u^k(x,t)\,dxdt  + \int_{\R^+} \beta^L(x)  u_0^k(x)\,dx, \label{halfline-fboundeq1}
\end{eqnarray}

\noindent from which, since $\uk \geq 0$,  it follows that 
\begin{equation}
\label{fpartbound}
k\int_0^T \int_1^L F(\uk, \vk)~dxdt
\end{equation}
\noindent is bounded independently of $L$, $k>0$,  since  (\ref{linfboundhalf}), the definition of $\beta^L$, and the fact that $\|u_0^k\|_{L^1(\R^+)}$ is bounded independently of $k$ imply that the right-hand side of \eqref{halfline-fboundeq1} is bounded independently of $k$ and $L$. On the other hand, multiplying the equation for $\vk$ by $\hat{\beta}$ and integrating over $\R^+ \times (0, t_0)$ yields
$$k \int_0^{t_0} \int_{\R^+} \hat {\beta} F(\uk, \vk)~dxdt = d_v \int_0^{t_0} \int_{\R^+} \hat{\beta}_{xx} \vk~dxdt - 
\int_{0}^2 \hat{\beta} [ \vk(x, T) - \vk_0(x)]~dx,$$
\noindent which, together with (\ref{linfboundhalf}), implies that
\begin{equation}
\label{fpartbound2}
k\int_0^T \int_0^1 F(\uk, \vk)~dxdt,
\end{equation}
is bounded independently of  $k>0$  and of $d_v \geq 0$ sufficiently small. The result then follows from (\ref{fpartbound}), (\ref{fpartbound2}), and Lebesgue's monotone convergence theorem. \qed

\bigskip
\begin{lemma}
\label{lemvliboundhalfline}
There exists a constant $C>0$,  independent of  $d_v \geq 0$ and $k>0$, such that for any  solution $(\uk, \vk)$ of $(P_2^k)$ satisfying (\ref{linfboundhalf}), we have
\begin{equation}
\label{vheatbound}
\int_0^{\infty} \uk(x, t_0) \\\, dx \; \leq \; C \;\; \mbox{and} \;\; \int_{0}^{\infty} |V_0 - \vk (x,t_0)|~dx \leq C\;\;\; \mbox{for all} \;\; t_0 \in [0, T].
\end{equation}
\end{lemma}

\pf The estimate for $\uk$ is immediate from \eqref{linfboundhalf}, \eqref{halfline-fboundeq1}, and Lebesgue's monotone convergence theorem.  
Then arguments similar to those used in the proof of Lemma \ref{lemloneboundwholeline} yield the estimate for $\vk$, since multiplying the equation satisfied by $\hat{v}^k:=V_0 - \vk$
by $m_{\alpha}'(\hat{v}^k) \psi^L_+$, where $\psi^L_+$ and $m_{\alpha}$ are as defined in the proofs of Lemmas \ref{lemcomparisonhalf} and \ref{lemloneboundwholeline}
respectively, integrating over $\R^+ \times (0, t_0)$, and noting that, if $d_v >0$, the boundary condition at $x=0$ yields
\begin{eqnarray*}
\int_0^{\infty} \psi^L_+ m_{\alpha}'(\hat{v}^k) \hat{v}^k_{xx} ~ dx & = & - \int_0^{\infty} (\psi^L_+)_x m_{\alpha}'(\hat{v}^k) \hat{v}^k_x + \psi^L_+ m_{\alpha}''(\hat{v}^k) (\hat{v}^k_x)^2~dx \\ & \leq & - \int_0^{\infty}( \psi^L_+)_x (m_{\alpha}(\hat{v}^k))_x~dx = \int_0^{\infty} (\psi^L_+)_{xx} m_{\alpha}(\hat{v}^k)~dx,
\end{eqnarray*}
\noindent together gives, after letting $\alpha \to 0$,  that for each $t_0 \in (0, T)$, 
\begin{eqnarray*}
\lefteqn{\int_0^{\infty} \psi^L_+ | \hat{v}^k(x, t_0)|dx \leq \int_0^{\infty} \psi^L_+ | \hat{v}^k(x,0)|dx \; } \\ & & \hspace*{3cm}+ \; d_v \int_0^{t_0} \int_0^{\infty} | \hat{v}^k| \left(\psi^L_+ \right)_{xx}dxdt + \int_0^{t_0}
  kF(\uk, V_0 - \hat{v}^k)\mbox{sgn}(\hat{v}^k) \psi^L_+ dxdt.
  \end{eqnarray*}
 \noindent The result then follows using (\ref{linfboundhalf}), Lemma \ref{lemfboundhalfline}, Lebesgue's monotone convergence theorem, and the fact that $\| V_0 - {v}^k_0\|_{L^1(\R^+)}$ is bounded independently of $k$. \qed

\bigskip
\noindent Next we prove the half-line analogue of Lemma \ref{leml2derivbound}. 

\begin{lemma}
\label{leml2bound2}
Suppose that $d_u >0$  and $d_v \geq 0$. Then  there exists $C>0$, independent of  $d_v$  and $k>0$, 
such that for   any solution $(u^k, v^k)$ of $(P_2^k)$
satisfying (\ref{linfboundhalf}),
\begin{equation}
\label{derivbounds2}
d_u \int_0^T \int_0^{\infty} (u^k_x)^2(x,t) ~dxdt \leq C \quad \mbox{and} \quad  d_v \int_0^T \int_{0}^{\infty} (v^k_x)^2(x,t)~dxdt \leq C.
\end{equation}
\end{lemma}

\pf Let $\hat{u} \in C^{\infty}(\R^+)$ be a fixed function such that $\hat{u}(0) =U_0$ and $\hat{u}(x) =0$ when $x \geq 1$. Define $y^k: = \uk - \hat{u}$.  Then
$y^k$ satisfies 
$$
 \left\{\begin{array}{rlrl}\displaystyle
y^k_t&= d_u y^k_{xx}  + d_u \hat{u}_{xx} - k F(y^k + \hat{u}, \vk) \quad &\mbox{in}~\st,
\\[.1 in]
\displaystyle y^k(0, t) &= 0 \quad &\mbox{for}~t \in [0, T], \\[.1 in]
 y^k(x,0)&=u^k_0(x) - \hat{u}(x),   &\qquad\mbox{for} \quad x\in \R^+,
\end{array}\right.$$
\noindent and multiplying the equation for $y^k$ by $y^k \psi^L_+$ and integrating over $S_T$ gives that 
\begin{eqnarray*}
\lefteqn{\frac{1}{2} \int_{0}^{\infty} \psi^L_+ (y^k)^2(x, T)dx + d_u \int_0^T  \int_{0}^{\infty}\psi^L_+ (y^k_x)^2dxdt
= \frac{1}{2} \int_{0}^{\infty} \psi^L_+ (u^k_0 -\hat{u})^2(x) dx}\\ &+ &  \frac{d_u}{2} \int_0^T  \int_{0}^{\infty} (y^k)^2 \left( \psi^L_+ \right)_{xx} dxdt
+ d_u \int_0^T \int_0^{\infty} \psi^L_+ \hat{u}_{xx} y^k dxdt
- k  \int_0^T  \int_{0}^{\infty} \psi^L_+ y^k F(\uk, \vk)dxdt.
\end{eqnarray*}
\noindent Since $y^k_x = \uk_x - \hat{u}_x$, the first estimate in (\ref{derivbounds2}) then follows using (\ref{linfboundhalf}), Lemma \ref{lemfboundhalfline} and the fact that $\|\uk_0\|_{L^1(\R^+)}$ is bounded independently of $k$. A similar argument yields the estimate for $\vk_x$, using the equation for $\hat{v}^k:=V_0 - \vk$ multiplied by $\psi^L_+ \hat{v}^k$
and the fact that $\hat{v}^k_x(0, t)=0$ for all $t \in (0, T)$. \qed

\bigskip

\noindent Recall the notation for space and time translates introduced in (\ref{spacetimedef}). 

\begin{lemma}
\label{lemspacetranslatehalf}
Suppose that $d_u > 0$ and $d_v \geq 0$, and let $(\uk, \vk)$  be a solution of $(P_2^k)$ satisfying (\ref{linfboundhalf}). Then for each \textcolor{black}{$r \in (0,1)$}, there exists a function $G_r \geq 0$, independent of $d_v \geq 0$ and $k >0$, such that $G_r(\xi) \to 0$ as $|\xi| \to 0$, and for all $|\xi| \leq \frac{r}{4}$ and $t \in (0, T)$, 
\begin{equation}
\label{spacetranslate}
\int_{r}^{\infty} | \uk(x,t)  - S_{\xi} \uk (x,t)| + | \vk (x,t)- S_{\xi} \vk(x,t)| ~dx \, \leq \, G_r(\xi).
\end{equation}
\end{lemma}

\pf Let $u,v,u_0$ and $v_0$ be as  defined in (\ref{ushiftdef}) and \textcolor{black}{define a cut-off function $\gamma^1_r \in C^{\infty}(\R^+)$ such that $ 0 \leq \gamma^1_r \leq 1$, $\gamma^1_r (x)=0$ when $x \in [0, \frac{r}{2}]$,  $\gamma^1_r (x)=1$ when $x \in [r, 1]$, and $\gamma^1_r (x) =0$ when $x \geq 2$. Then given $L \geq 1$, define the family of cut-off functions
 $\gamma^L_r \in C^{\infty}(\R^+)$ by $\gamma^L_r(x) = \gamma^1_r(x)$ when $x \in [0, r]$, $\gamma^L_r(x) = 1$ when $x \in [r, L]$, and $\gamma^L_r(x) = \gamma^1_r(x+1-L)$ when $x \geq L$.
 Note that $0 \leq \gamma^L_r \leq 1$ for all $L$, and 
 $\left(\gamma^L_r\right)_x$, $\left(\gamma^L_r\right)_{xx}$ are bounded in both $L^{\infty}(\R^+)$ and $L^1(\R^+)$ independently of $L$. }
Then 
$$
\begin{array}{rlrl}\displaystyle
u_t&= d_u u_{xx} -k\{ F(\uk, \vk) - F(S_{\xi}\uk, S_{\xi}\vk)\} \quad &\mbox{in}~(\textstyle{\frac{r}{4}}, \infty) \times (0, T), 
\\[.1 in]
\displaystyle v_t&= d_v v_{xx} -k\{ F(\uk, \vk) - F(S_{\xi}\uk, S_{\xi}\vk)\} \quad &\mbox{in}~(\textstyle{\frac{r}{4}}, \infty) \times (0, T), \\[.1 in]
 u(x,0)&= \uk_0(x) - \uk_0(x + \xi) ,\quad v(x,0)=\vk_0(x) - \vk_0(x + \xi)   &\qquad\mbox{for} \quad x\in (\textstyle{\frac{r}{4}}, \infty),
\end{array}$$

\noindent so that arguing as in the proof of Lemma \ref{leml1contraction} yields that for each $t_0 \in (0, T)$, 
\begin{eqnarray}
\label{jumptoit}
\lefteqn{\;\;\int_{\frac{r}{2}}^{\infty} \gamma^L_r (x) \{ |u(x, t_0)| + |v(x, t_0)| \}~dx \leq \int_{\frac{r}{2}}^{\infty} \gamma^L_r (x) \{ |u_0(x)| + | v_0(x)|\}~dx} \\ & & \hspace*{5.5cm} + \int_0^{t_0} \int_{\frac{r}{2}}^{\infty} \left(\gamma^L_r\right)_{xx}(x) \{ d_u|u(x,t)| + d_v |v(x,t)| \} ~dxdt.  \nonumber
\end{eqnarray}

\noindent Now, by the definition of $u$, 
\begin{eqnarray*}
\int_{\frac{r}{2}}^{\infty} \left(\gamma^L_r\right)_{xx}(x)  d_u |u(x,t)| dx & = & \int_{\frac{r}{2}}^{\infty} \left(\gamma^L_r\right)_{xx}(x)  d_u | 
\uk(x,t) - \uk(x+\xi, t)| dx \\
& = & \int_{\frac{r}{2}}^{\infty} \left(\gamma^L_r\right)_{xx}(x)  d_u \left| \int_0^1 \uk_x(x+ \theta \xi, t) \xi d \theta  \right| dx \\ & \leq & |\xi| d_u 
\int_{\frac{r}{2}}^{\infty} \left(\gamma^L_r\right)_{xx}(x)  \int_0^1 |\uk_x(x+ \theta \xi, t) | d \theta dx,
\end{eqnarray*}

\noindent thus
\begin{eqnarray*}
\lefteqn{\int_0^{t_0} \int_{\frac{r}{2}}^{\infty} \left(\gamma^L_r\right)_{xx}(x)  d_u |u(x,t)| dxdt } \\ & & \leq  
|\xi| \sqrt{d_u} 
\int_0^1 \left( \int_0^{t_0} \int_{\frac{r}{2}}^{\infty} \left(\gamma^L_r\right)_{xx}^2 dxdt  \right)^{\frac{1}{2}}
\left( \int_0^{t_0} \int_{\frac{r}{2}}^{\infty} d_u |\uk_x(x+ \theta \xi, t) |^2  dxdt  \right)^{\frac{1}{2}} d \theta \\
& & \leq |\xi| \sqrt{d_u} \left( \int_0^{t_0} \int_{\frac{r}{2}}^{\infty} \left(\gamma^L_r\right)_{xx}^2 dxdt  \right)^{\frac{1}{2}}
\left( \int_0^{t_0} \int_{\frac{r}{4}}^{\infty} d_u |\uk_x(x, t) |^2  dxdt  \right)^{\frac{1}{2}},
\end{eqnarray*}

\noindent and hence applying Lemma \ref{leml2bound2}  shows that 
\begin{equation}
\label{estr}
\int_0^{t_0} \int_{\frac{r}{2}}^{\infty} \left(\gamma^L_r\right)_{xx}(x)  d_u |u(x,t)| dxdt  \leq K_r |\xi|,
\end{equation}
\noindent for some constant $K_r$. The result then follows from (\ref{jumptoit}) using (\ref{estr}),  a similar estimate for $\int_0^{t_0} \int_{\frac{r}{2}}^{\infty} \left(\gamma^L_r\right)_{xx}  d_v |v| dxdt $,
the fact that $\| \uk_0 (\cdot + \xi) - \uk_0 (\cdot) \|_{L^1((r, \infty))} + \| \vk_0 (\cdot + \xi) - \vk_0 (\cdot) \|_{L^1((r, \infty))} \; \leq \; \omega_r(|\xi|)$ where $\omega_r(|\xi|) \to 0$ as $\xi \to 0$,
 and Lebesgue's monotone convergence theorem. \qed

\bigskip
\begin{lemma}
\label{lemtimetranslateshalf}
Suppose that $d_u > 0$, $d_v \geq 0$ and let $(u^k, v^k)$ be a solution of $(P_2^k)$
satisfying (\ref{linfboundhalf}). Then there exists $C>0$, independent of  $d_v$ and $k$, such
that for any $\tau \in (0, T)$, 
\begin{eqnarray*}
\int_0^{T- \tau} \int_{0}^{\infty} | T_{\tau}\uk(x,t)- u^k(x,t)|^2 ~ dx dt & \leq & \tau  C,\\
\int_0^{T- \tau} \int_{0}^{\infty} | T_{\tau}\vk(x,t) - v^k(x,t)|^2 ~ dx dt & \leq & \tau C.
\end{eqnarray*}
\end{lemma}

\pf This follows from arguments analogous to those used in the proof of Lemma \ref{lemtimetranslates}, replacing  $\psi^L$ by $\psi^L_+ := \psi^L|_{\R^+}$ and integrals over $\R$ by integrals over $\R^+$, noting that $\uk(0, t+ \tau) - \uk(0, t)=0$  and $\vk_x(0, t+ \tau) - \vk_x(0, t)=0$ for all $t \in (0, T- \tau)$  and using the bounds in Lemmas 
\ref{lemfboundhalfline} and \ref{leml2bound2}. \qed

\bigskip

\begin{lemma}
\label{lemconvergencedzerohalf} Let $k>0$ and $d_u>0$ be fixed and $(u^k_{d_v}, v^k_{d_v})$ be solutions of $(P_2^k)$ satisfying (\ref{linfboundhalf}) with $d_v >0$. Then there exists $(u^k_*, v^k_*) \in (L^{\infty}(S_T))^2$ such that up to a subsequence, for each $J>0$, 
$$ \begin{array}{cl}
u^k_{d_v} \to u^k_* & \mbox{in} \;\; L^2((0, J) \times (0, T)),\\
v^k_{d_v} \to v^k_* & \mbox{in} \;\; L^2((0, J) \times (0, T)),\\
u^k_{d_v} - \hat{u} \rightharpoonup u^k_* - \hat{u} & \mbox{in} \;\; L^2(0, T; H^1_0(\R^+)),
\end{array}
$$
\noindent as $d_v \to 0$, where $\hat{u} \in C^{\infty}(\R^+)$ is a smooth function such that $\hat{u}(0) = U_0$ and $\hat{u}(x) = u_0^{\infty}(x)$ for all $x \geq 1$. 
\end{lemma}

\pf It follows from (\ref{linfboundhalf}) and Lemma \ref{lemvliboundhalfline}  that $\|u^k_{d_v} - u^{\infty}_0\|_{L^2(S_T)}$ and $\|v^k_{d_v} - v^{\infty}_0\|_{L^2(S_T)}$ are bounded independently of $d_v$. So  Lemmas \ref{lemspacetranslatehalf}, \ref{lemtimetranslateshalf} and the Riesz-Fr\'echet-Kolmogorov Theorem \cite[Theorem 4.26]{brezisnew} yield that the sets $\{ u^k_{d_v} - u^{\infty}_0\}_{d_v>0}$ and
$\{ v^k_{d_v} - v^{\infty}_0\}_{d_v>0}$ are each relatively compact in $L^2((0, J) \times (0, T))$ for each $J>0$. The weak convergence of $u^k_{d_v} - \hat{u}$ in $L^2(0, T; H^1_0(\R))$ follows from  the fact that $\|u^k_{d_v} \|_{L^2(S_T)}$ is bounded independently of $d_v$ together with the proof of Lemma \ref{leml2bound2}. \qed

\bigskip

\noindent Lemma \ref{lemconvergencedzerohalf}  and Corollary \ref{cor-half} enable the following result to be established using arguments similar to those that yield Theorem \ref{thmexistzero}. We omit details
of the proof. 

\begin{theorem}
\label{thmexistzerohalf}
Let $d_v=0$ and $k>0$. Then problem $(P_2^k)$ has a unique weak solution
$$(\uk, \vk) \in W^{2,1}_p( (0, J) \times (0, T)) \times W^{1,\infty}(0, T; L^{\infty}((0, J)))\;\; \mbox{for each} \; J>0 , \; \; p \geq 1,$$
where $(\uk, \vk)$ is a weak solution in the sense that 
\begin{eqnarray}
\;\;\iint_{S_T} \uk \psi_t ~dxdt  + \iint_{S_T} \{ d_u \uk \psi_{xx} - k F(\uk, \vk) \psi \}~dxdt & =  & - \int_{0}^{\infty} u_0^k \psi( \cdot, 0) ~ dx   \label{udzhalf} \\
& &  \hspace*{1cm}   -  d_u U_0 \int_0^T \psi_x (0,t)~dt 
, \nonumber  \\
\iint_{Q_T} \vk \psi_t ~dxdt   -  \iint_{Q_T} k F(\uk, \vk) \psi ~dxdt  & =  & - \int_{\R} v_0^k
\psi( \cdot, 0) ~ dx, \label{vdzhalf}
\end{eqnarray}
for all $\psi \in \hat{\mathcal{F}}_T = \{\psi \in C^{2,1}(S_T): \psi (\cdot, T) =0,\; \psi(0, t) =0 \; \mbox{for} \; t \in (0, T)\; \mbox{and} \;  \mbox{supp}\, \psi \subset [0, J] \times [0, T] \;\mbox{for some} \; J>0 \}$, and also satisfies  $0 \leq \uk, \vk \leq M$.

\end{theorem}

\label{subsectiondvzerohalf}

\subsection{The limit problem for $(P_2^k)$ as $k \to \infty$}

\noindent The next result follows directly from arguments similar to those used in section \ref{subsectionkinftywhole}, exploiting the half-line estimates established in section \ref{subsectiondvzerohalf}. 

\begin{lemma}
\label{lemkinftyconvergencehalf}
Let $d_u >0$ and $d_v \geq 0$ be fixed and $(u^k, v^k)$ be solutions of $(P_2^k)$ satisfying (\ref{linfboundhalf}) with $k>0$. Then there exists $(u, v) \in (L^{\infty}(S_T))^2$ such that up to a subsequence, for each $J>0$, 
$$ \begin{array}{cl}
u^k \to u & \mbox{in} \;\; L^2((0, J) \times (0, T)),\\
v^k \to v & \mbox{in} \;\; L^2((0, J) \times (0, T)),\\
u^k - \hat{u} \rightharpoonup u - \hat{u} & \mbox{in} \;\; L^2(0, T; H^1_0(\R)),
\end{array}
$$
\noindent as $k \to \infty$, where $\hat{u} \in C^{\infty}(\R^+)$ is a smooth function such that $\hat{u}(0) = U_0$ and $\hat{u}(x) =0$ for all $x \geq 1$. Moreover, 
\begin{equation}
\label{seghalf}
uv =0 \;\;\; a.e. \; \mbox{in} \;\; S_T.
\end{equation}

\end{lemma}

\medskip
\noindent Taking $\wk$ and $w$ as in (\ref{wkwdef}), we clearly again have that as a sequence $k_n \to \infty$, 
$w^{k_n} \to w$ in $L^2(S_T)$ and almost everywhere in $S_T$, 
and that $u=w^+ $ and $v = -w^-$.  

\medskip \noindent The next result is our half-line counterpart of Lemma \ref{lemlimiteq}, with the boundary at $x=0$ clearly now playing a r\^{o}le. Note that  here, similarly to \cite{hmm},  the limit function $w$ satisfies a Dirichlet boundary condition  both  when $d_v=0$ and when $d_v>0$.

\begin{lemma}
\label{lemlimiteqhalf}
Let $d_u >0$, $d_v \geq 0$ and $(u,v)$ be as in Lemma \ref{lemkinftyconvergencehalf}. Then
\begin{equation}
\label{limiteqhalf}
- \iint_{S_T} (u-v) \psi_t~dxdt  - \int_{0}^{\infty} (u^{\infty}_0 - v^{\infty}_0) \, \psi(x, 0)~dx = d_u U_0 \int_0^T \psi_x(0, t) dt + \iint_{S_T}(d_u u - d_v v) \
\psi_{xx} ~dxdt,
\end{equation}
for all $\psi \in \hat{\mathcal{F}}_T = \{\psi \in C^{2,1}(S_T): \psi (\cdot, T) =0,\; \psi(0, t) =0 \; \mbox{for} \; t \in (0, T)\; \mbox{and} \;  \mbox{supp}\, \psi \subset [0, J] \times [0, T] \;\mbox{for some} \; J>0 \}$. 

\end{lemma}

\pf Multiplying the difference between the equations for $\uk$ and $\vk$ by $\psi \in \hat{\mathcal{F}}_T$ and integrating over $S_T$ gives
\begin{eqnarray}
\lefteqn{- \iint_{S_T} (\uk-\vk) \psi_t~dxdt  - \int_{0}^{\infty} (u^k_0 - v^k_0) \, \psi(x, 0)~dx}\label{firsthalf} \\ & &  \hspace*{1cm} =  \int_0^T \{d_u \uk(0, t)  - d_v \vk(0, t)  \}\psi_x(0, t) dt  + \iint_{S_T}(d_u \uk - d_v \vk) 
\psi_{xx} ~dxdt. \nonumber
\end{eqnarray}

\noindent Now it follows exactly as argued in the proof of \cite[Prop. 8]{hmm} that if $d_v>0$, then the segregation property (\ref{seghalf}) yields that as $k \rightarrow \infty$, 
$$\gamma(d_u\uk - d_v \vk) \rightharpoonup d_u U_0 \;\;\; \mbox{in} \;\; L^2(\{0\} \times (0, T)),$$
\noindent where $\gamma$ denotes the trace on the boundary $\{0\} \times (0, T)$. So (\ref{limiteqhalf}) follows by letting $k \to \infty$ in (\ref{firsthalf}). \qed

\bigskip
\noindent Now recall the definition of $\mathcal{D}$ from (\ref{defcurld}) and define
 the limit problem
 \medskip
$$(P_2^{limit}) \left\{\begin{array}{rlrl}\displaystyle
w_t&= \mathcal{D}(w)_{xx},\quad &\mbox{in}~[0, \infty) \times (0, \infty),
 \\[.1 in]
 w(x,0)&= w_0(x) : =  - V_0, & \mbox{if} \;\;x <0, \\
 w(0,t) &= U_0,& \mbox{for}\; \;t \in (0, \infty).   \end{array} \right.$$
 
 \bigskip
 
 \begin{definition}
 \label{defweaksolutionhalf}
 A function $w$ is a weak solution of Problem $(P_2^{limit})$ if
 \begin{itemize}
 \item[(i)] $w \in L^{\infty} (\R^+ \times \R^+)$,
 \item[(ii)] for all $T>0$, 
 $$\iint_{S_T} (w \psi_t + \mathcal{D}(w) \psi_{xx})~dxdt = - d_u U_0 \int_0^T \psi_x(0, t) dt - \int_{0}^{\infty} w_0(x) \psi(x,0)~dx,$$
 for all $\psi \in \hat{\mathcal{F}}_T = \{\psi \in C^{2,1}(S_T): \psi (\cdot, T) =0,\; \psi(0, t) =0 \; \mbox{for} \; t \in (0, T)\; \mbox{and} \;  \mbox{supp}\, \psi \subset [0, J] \times [0, T] \;\mbox{for some} \; J>0 \}$. 
  \end{itemize}
 \end{definition}
 
 \bigskip
 
 \begin{lemma}
 \label{lemuniqueweaklimithalf}
 With $(u,v)$ from Lemma \ref{lemkinftyconvergencehalf},  the function $w:u-v$ is the unique weak solution of Problem $(P_2^{limit})$ and the whole sequence $(\uk, \vk)$ in Lemma
 \ref{lemkinftyconvergencehalf} converges to $(w^+, -w^-)$.
 \end{lemma}
 
 \pf That $w$ is a weak solution of $(P_2^{limit})$ follows immediately from Lemma \ref{lemlimiteqhalf} and the definition of $\mathcal{D}$. 
 Minor modifications in the arguments used to establish \cite[Appendix, Proposition A]{bkp} and \cite[Proposition 9]{acp} yield a corresponding estimate with the domain $\R$ in \cite[Appendix, Proposition A]{bkp}  replaced by $\R^+$, from which uniqueness again follows via the reasoning in the proof of \cite[Appendix, Proof of Theorem C]{bkp}. \qed
 
 \bigskip
 \noindent As in the whole-line case, we can identify the limit $w$ as a certain self-similar solution both when $d_v>0$ and when $d_v=0$. We first state the analogue of  Corollary \ref{corfreeboundary}. 
 
% \begin{theorem}
%\label{thmfreeboundaryhalf}
%Let $w$ be the unique weak solution of Problem $(P_2^{limit})$. Suppose that there exists a function $\xi: [0, T] \to \R^+$ such that  for each $t \in [0, T]$, 
%$$ w(x, t) > 0 \;\;\; \mbox{if} \; \;\; x < \xi(t) \;\;\; \mbox{and} \;\;\; w(x, t) < 0 \;\;\; \mbox{if} \;\;\; x> \xi(t). $$
%Then if $t \mapsto \xi(t)$ is sufficiently smooth and the functions $u := w^+$ and $v: = - w^-$ are smooth up to $\xi(t)$, the functions
%$u$ and $v$ satisfy 
%$$ 
%(P_2^{limit}) \left\{
%\begin{array}{ll}
%u_t = d_u u_{xx},  & \mbox{in} \;\; \{(x,t) \in Q_T: x <  \xi(t)\} ,\\
%v_t = d_v v_{xx},   & \mbox{in} \;\; \{(x,t) \in Q_T: x >  \xi(t)\} ,\\
%\left[u \right] = d_v \left[v \right] =0,     & \mbox{on} \;\; \Gamma_T  := \{(x,t) \in Q_T: x =  \xi(t)\},\\
%\left [v \right] \xi'(t) = \left[ d_u u_x - d_v v_x \right], & \mbox{on} \;\; \Gamma_T := \{(x,t) \in Q_T: x =  \xi(t)\},\\
%u=U_0, & \mbox{on} \;\; \{0\} \times [0, T],\\
%u(\cdot, 0) = 0, & \mbox{in} \;\; (0, \infty),\\
%v(\cdot, 0) = V_0,  & \mbox{in} \;\; (0, \infty),
%\end{array}
%\right.
%$$

%\noindent where $\left[ \cdot \right]$ denotes the jump across $\xi(t)$ from $\{x < \xi(t) \}$ to $\{x > \xi(t)\}$, $\xi'(t)$ denotes the speed of propagation of the free boundary $\xi (t)$, and  we adopt
%the notational convention that a term multiplied by $d_v$ is considered to be absent if $d_v=0$.
%\end{theorem}

\bigskip

\begin{proposition}
\label{propfreeboundaryhalf}
Let $w$ be the unique weak solution of Problem $(P_2^{limit})$. Suppose that there exists a function $\xi: [0, T] \to \R^+$ such that  for each $t \in [0, T]$, 
$$ w(x, t) > 0 \;\;\; \mbox{if} \; \;\; x < \xi(t) \;\;\; \mbox{and} \;\;\; w(x, t) < 0 \;\;\; \mbox{if} \;\;\; x> \xi(t). $$
Then if $t \mapsto \xi(t)$ is sufficiently smooth and the functions $u := w^+$ and $v: = - w^-$ are smooth up to $\xi(t)$, the functions
$u$ and $v$ 
 satisfy one of two limit problems, depending on whether $d_v>0$ or $d_v=0$. If $d_v>0$, then 
$$ 
(P^{limit}_{2, d_v>0}) \left\{
\begin{array}{ll}
u_t = d_u u_{xx},  & \mbox{in} \;\; \{(x,t) \in Q_T: x <  \xi(t)\} ,\\
v = 0,  & \mbox{in} \;\; \{(x,t) \in Q_T: x <  \xi(t)\} ,\\
v_t = d_v v_{xx},   & \mbox{in} \;\; \{(x,t) \in Q_T: x >  \xi(t)\} ,\\
u = 0,   & \mbox{in} \;\; \{(x,t) \in Q_T: x >  \xi(t)\} , \\
\lim_{x \uparrow \xi(t)} u(x,t) =0 = \lim_{x \downarrow \xi(t)} v(x,t),    & \mbox{for each} \;\; t \in [0, T], \\
d_u \lim_{x \uparrow \xi(t)} u_x(x,t) = - d_v \lim_{x \downarrow \xi(t)} v_x(x,t), & \mbox{for each} \;\; t \in [0, T],\\
u=U_0, & \mbox{on} \;\; \{0\} \times [0, T],\\
u(\cdot, 0) = 0, & \mbox{in} \;\; (0, \infty),\\
v(\cdot, 0) = V_0,  & \mbox{in} \;\; (0, \infty),
\end{array}
\right.
$$
\noindent whereas if $d_v=0$ and  we suppose additionally that 
$\xi(0)=0$ and $t \mapsto \xi(t)$ is a non-decreasing function, then 
$$ 
(P^{limit}_{2, d_v=0}) \left\{
\begin{array}{ll}
u_t = d_u u_{xx},  & \mbox{in} \;\; \{(x,t) \in Q_T: x <  \xi(t)\} ,\\
v = 0,  & \mbox{in} \;\; \{(x,t) \in Q_T: x <  \xi(t)\} ,\\
v = V_0,   & \mbox{in} \;\; \{(x,t) \in Q_T: x >  \xi(t)\} ,\\
u = 0,   & \mbox{in} \;\; \{(x,t) \in Q_T: x >  \xi(t)\} ,\\
\lim_{x \uparrow \xi(t)} u(x,t) =0,  & \mbox{for each} \;\; t \in [0, T],   \\
V_0 \; \xi'(t) = - {d_u}\lim_{x \uparrow \xi(t)} u_x (x,t), & \mbox{for each} \;\; t \in [0, T],\\
u=U_0, & \mbox{on} \;\; \{0\} \times [0, T],\\
u(\cdot, 0) = 0, & \mbox{in} \;\; (0, \infty),\\
v(\cdot, 0) = V_0,  & \mbox{in} \;\; (0, \infty),
\end{array}
\right.
$$
where 
$\xi'(t)$ denotes the speed of propagation of the free boundary $\xi (t)$. 
\end{proposition}

\bigskip

\begin{theorem}
\label{thmselfsimilarhalf}
The unique weak solution $w$ of Problem $(P_2^{limit})$ has a self-similar form. There exists a function $f: \R^+ \to \R$ and a positive constant $a>0$ such that
\begin{equation}
\label{wform2}
w(x,t) = f\left( \frac{x}{\sqrt{t}}\right) , \;\;\; (x,t) \in S_T, \;\;\; \mbox{and} \;\;\;\xi(t) = a \sqrt{t}, \;\;\; t \in [0, T].
\end{equation}
\noindent If $d_v>0$, then $a >0 $ is the unique root of the equation
$$  d_u U_0 \int_a^{\infty} e^{\frac{a^2 - s^2}{4d_v}} ds = d_v V_0 \int_{0}^{a} e^{\frac{a^2 - s^2}{4d_u}},$$
\noindent and
\begin{equation}
\label{fonehalf}
f(\eta) = \left\{ \begin{array}{ll} U_0 \left(  1 - \frac{\int_{0}^{\eta} e^{- \frac{s^2}{4d_u}}\;ds}{\int_{0}^{a} e^{- \frac{s^2}{4d_u}} \;ds}        \right), & \mbox{if} \;\; \eta \leq a,\\
-V_0 \left(  1 - \frac{\int_{\eta}^{\infty} e^{- \frac{s^2}{4d_v}}\;ds}{\int_{a}^{\infty} e^{- \frac{s^2}{4d_v}} \;ds}        \right), & \mbox{if} \;\; \eta>a.
 \end{array} \right.
 \end{equation}
\noindent On the other hand, if $d_v=0$, then $a>0$ is the unique root of the equation
$$U_0 = \frac{V_0 a}{2d_u} \int_{0}^a e^{\frac{a^2 - s^2}{4 d_u}}\;ds,$$
\noindent and
\begin{equation}
\label{ftwohalf}
f(\eta) = \left\{ \begin{array}{ll} U_0 \left(  1 - \frac{\int_{0}^{\eta} e^{- \frac{s^2}{4d_u}}\;ds}{\int_{0}^{a} e^{- \frac{s^2}{4d_u}} \;ds}        \right), & \mbox{if} \;\; \eta \leq a,\\
-V_0 , & \mbox{if} \;\; \eta>a.
 \end{array} \right.
 \end{equation}
\end{theorem}

\label{seckliminthalf}

\bigskip
\section{Long-time behaviour for $(P_1^k)$ and $(P_2^k)$ when $k$ is fixed}

\noindent We conclude our study  by exploiting a scaling argument, due first to Kamin \cite{Kamin} and used also in \cite{hhp1, hhp2}, to infer the self-similar $t \to \infty$  limits of solutions $(\uk, \vk)$ of $(P_1^k)$ or $(P_2^k)$  from the $k \to \infty$ limits discussed in Sections \ref{sec-klimitwhole} and \ref{seckliminthalf}. As mentioned in the Introduction, this enables us, in particular, to give rigorous justification to the long-time asymptotics of reaction fronts discussed by Trevelyan et al \cite{TSDW}.
Note that in \cite{hhp1, hhp2}, uniform convergence results  as $k \to \infty$ implied pointwise convergence results as $ t \to \infty$ for solutions of a problem with $k$ fixed. Here, however,  we simply use the $L^2$-convergence from Lemmas 
\ref{lemkinftyconvergence} and \ref{lemkinftyconvergencehalf} to deduce convergence in a certain average sense of solutions of $(P_1^k)$ and $(P_2^k)$ as $ t \to \infty$ along a subsequence.

\begin{theorem}
\label{thmtinftywhole}
Let $(\uk, \vk)$ be the solution of problem $(P_1^k)$ with initial data $\uk_0$, $\vk_0 \in C^2(\R)$ such that
\begin{equation}
\label{ic1}
\|\uk_0 - u^{\infty}_0\|_{L^1(\R)} < \infty, \;\;\; \|\vk_0 - v^{\infty}_0\|_{L^1(\R)} < \infty,
\end{equation}
and
\begin{equation}
\label{ic2}
\uk_0(x) \rightarrow U_0, 0 \;\; \mbox{as} \;\;x \rightarrow -\infty, \infty \;\; \mbox{and} \;\; 
\vk_0(x) \rightarrow 0, V_0 \;\; \mbox{as} \;\;x \rightarrow -\infty, \infty.
\end{equation}

\bigskip
\noindent Then for each $J>0$, there exists a sequence $t_n \to \infty$ such that
\begin{equation}
\label{longtimeconvu}
\frac{1}{\sqrt{t_n}} \int_{- J \sqrt{t_n}}^{J \sqrt{t_n}} \left| \uk(y, t_n) - f^+\left( \frac{y}{\sqrt{t_n}} \right)   \right|^2 ~ dy \;\; \to \;\; 0 \;\; \mbox{as} \;\; t_n \to \infty,
\end{equation}
and
\begin{equation}
\label{longtimeconvv}
\frac{1}{\sqrt{t_n}} \int_{- J \sqrt{t_n}}^{J \sqrt{t_n}} \left| \vk(y, t_n) + f^-\left( \frac{y}{\sqrt{t_n}} \right)   \right|^2 ~ dy \;\; \to \;\; 0 \;\; \mbox{as} \;\; t_n \to \infty.
\end{equation}

\bigskip
\noindent where $f$ is the self-similar profile given by (\ref{fone}) if $d_v>0$, and by (\ref{ftwo}) if $d_v=0$. Here,  as usual, 
$f^+ := \max \{f, 0\}$, $f^- := \min\{f, 0\}$. 
\end{theorem}

\bigskip
\pf For each $l>0$, the scaled functions 
$$\uk_l (x,t): =  \uk(lx, l^2 t), \;\;\; \vk_l(x,t) := \vk_l(x,t) = \vk(lx, l^2t),$$
satisfy the system
$$
(P1^k_l) \left\{\begin{array}{rlrl}\displaystyle
u_t&= d_u u_{xx} -kl^2F(u,v)\quad &\mbox{in}~Q_{l^2T},
\\[.1 in]
\displaystyle v_t&= d_v v_{xx} -kl^2F(u,v)\quad &\mbox{in}~Q_{l^2T}, \\[.1 in]
 u(x,0)&=u^k_0(lx),\quad v(x,0)=v^k_0(lx)   &\qquad\mbox{for} \quad x\in \R.
\end{array}\right.$$
Moreover, it follows from (\ref{ic1}) that as $l \to \infty$, 
$$\|\uk_l(\cdot, 0) - u_0^{\infty}\|_{L^1(\R)} \to 0, \;\;\; \|\vk_l(\cdot, 0) - v_0^{\infty}\|_{L^1(\R)} \to 0.$$
So Lemma \ref{lemkinftyconvergence}, Lemma \ref{lemuniqueweaklimit} and 
Theorem \ref{thmselfsimilar} together imply that for each $T>0$ and each $J>0$, 
\begin{equation}
\label{basiclconv}
\int_0^T \int_{-J}^J \left| \uk_l(x,t) - f^+\left( \frac{x}{\sqrt{t}} \right) \right|^2~dxdt \to 0, \;\;\;
\int_0^T \int_{-J}^J \left| \vk_l(x,t) + f^-\left( \frac{x}{\sqrt{t}} \right) \right|^2~dxdt \to 0 \;\; \mbox{as} \;\;
l \to \infty,
\end{equation}

\noindent where $f$ is as in the statement of the theorem. Now take $T \geq 1$ and fix $J>0$. Then it follows from
(\ref{basiclconv}) that there exists $t_0 \in (\frac{1}{2}, 1)$ and a sequence $l_n \to \infty$ such that as $l_n \to \infty$, 
$$\int_{-J}^J \left| \uk_{l_n}(x,t_0) - f^+\left( \frac{x}{\sqrt{t_0}} \right) \right|^2~dx \to 0, \;\;\;
\int_{-J}^J \left| \vk_{l_n}(x,t_0) + f^-\left( \frac{x}{\sqrt{t_0}} \right) \right|^2~dx \to 0,$$
\noindent or equivalently
$$\int_{-J}^J \left| \uk(l_nx,l_n^2t_0) - f^+\left( \frac{x}{\sqrt{t_0}} \right) \right|^2~dx \to 0, \;\;\;
\int_{-J}^J \left| \vk(l_n x,l_n^2t_0) + f^-\left( \frac{x}{\sqrt{t_0}} \right) \right|^2~dx \to 0,$$
\noindent which yields immediately that as $ l_n \to \infty$, 
\begin{equation}
\label{penum}
\frac{1}{l_n}\int_{-l_nJ}^{l_nJ} \left| \uk(y,l_n^2t_0) - f^+\left( \frac{y}{l_n \sqrt{t_0}} \right) \right|^2~dy \to 0, \;\;\;
\frac{1}{l_n}\int_{-l_nJ}^{l_nJ} \left| \vk(y,l_n^2t_0) + f^-\left( \frac{y}{l_n \sqrt{t_0}} \right) \right|^2~dy \to 0,
\end{equation}
\noindent Taking $s_n:= l_n^2 t_0$ in (\ref{penum}) then gives that
$$\sqrt{\frac{t_0}{s_n}}\int_{-J\sqrt{\frac{s_n}{t_0}}}^{J\sqrt{\frac{s_n}{t_0}}} \left| \uk(y,s_n) - f^+\left( \frac{y}{\sqrt{s_n}} \right) \right|^2~dy \to 0, \;\;\;
\sqrt{\frac{t_0}{s_n}}\int_{-J\sqrt{\frac{s_n}{t_0}}}^{J\sqrt{\frac{s_n}{t_0}}} \left| \vk(y,s_n) + f^-\left( \frac{y}{\sqrt{s_n}} \right) \right|^2~dy \to 0,$$
\noindent as $s_n \to \infty$, from which the result follows. \qed

\bigskip
\noindent Minor modifications in the  arguments above show the following, corresponding result for the half-line problem $(P_2^k)$. 
We leave the details of the proof to the reader. 

\begin{theorem}
\label{thmtinftyhalf}
Let $(\uk, \vk)$ be the solution of problem $(P_2^k)$ with initial data $\uk_0$, $\vk_0 \in C^2(\R^+)$ such that
\begin{equation}
\label{ic1half}
\|\uk_0 - u^{\infty}_0\|_{L^1(\R^+)} < \infty, \;\;\; \|\vk_0 - v^{\infty}_0\|_{L^1(\R^+)} < \infty,
\end{equation}
and
\begin{equation}
\label{ic2half}
\uk_0(x) \rightarrow  0 \;\; \mbox{as} \;\;x \rightarrow \infty \;\; \mbox{and} \;\; 
\vk_0(x) \rightarrow  V_0 \;\; \mbox{as} \;\;x \rightarrow \infty.
\end{equation}

\bigskip
\noindent Then for each $J>0$, there exists a sequence $t_n \to \infty$ such that
\begin{equation}
\label{longtimeconvuhalf}
\frac{1}{\sqrt{t_n}} \int_{0}^{J \sqrt{t_n}} \left| \uk(y, t_n) - f^+\left( \frac{y}{\sqrt{t_n}} \right)   \right|^2 ~ dy \;\; \to \;\; 0 \;\; \mbox{as} \;\; t_n \to \infty,
\end{equation}
and
\begin{equation}
\label{longtimeconvvhalf}
\frac{1}{\sqrt{t_n}} \int_{0 }^{J \sqrt{t_n}} \left| \vk(y, t_n) + f^-\left( \frac{y}{\sqrt{t_n}} \right)   \right|^2 ~ dy \;\; \to \;\; 0 \;\; \mbox{as} \;\; t_n \to \infty.
\end{equation}

\bigskip
\noindent where $f$ is the self-similar profile given by (\ref{fonehalf}) if $d_v>0$, and by (\ref{ftwohalf}) if $d_v=0$. \end{theorem}

\bigskip
\noindent {\bf Acknowledgements}
We would like to thank Philip Trevelyan, whose asymptotic work on reaction fronts was one of the motivations of our research, for discussions on the physical background to the systems studied, \textcolor{black}{and the referee for his/her useful comments}. We are  grateful also to  Fran\c{c}ois Hamel, Shoshana Kamin, Hiroshi Matano and Dariusz Wrzosek for interesting remarks on our results, which have improved our exposition and suggested potential directions for further developments. This work was partly supported by a grant from the London Mathematical Society.

%The function $\hat{v}^k:= V_0 - \vk$ satisfies 
%$$
 %\left\{\begin{array}{rlrl}\displaystyle
%\hat{v}^k_t&= d_v \hat{v}^k_{xx} - k F(\uk, V_0 - \hat{v}^k), \quad &\mbox{in}~\st,
%\\[.1 in]
%\displaystyle \hat{v}^k_x(0, t) &= 0 \quad &\mbox{for all}~t \in [0, T], \\[.1 in]
 %\hat{v}^k (x,0)&=V_0 - v^k_0(x),   &\qquad\mbox{for} \quad x\in \R^+.
%\end{array}\right.$$
%\noindent Then taking $\psi^L_+$ as in the proof of Lemma \ref{lemcomparisonhalf}

\end{document}